\documentclass[12 pt]{article}
\usepackage{amsmath,amsfonts,amssymb,amsthm}
\usepackage{latexsym}
\usepackage[dvips]{epsfig}
\usepackage{amscd}
\usepackage{pstricks,pst-node}

\def\C{\mathbb{C}}
\def\Z{\mathbb{Z}}
\def\N{\mathbb{N}}

\def\R{{\mathbf{R}}}
\def\RR{\mathbb{R}}
\def\qed{$\hfill \square$}
\def\G{\Gamma}
\def\g{\mathfrak{g}}
\def\n{\mathfrak{n}}
\def\h{\mathfrak{h}}
\def\V{\mathbb{V}}
\def\W{\mathbb{W}}
\def\R{\mathbb R}

\def\sotimes{\stackrel{\cdot}{\otimes}}

\def\lra{\longrightarrow}

\def\Mod{\mathbf{Mod}}
\def\UMod{\underline{\mathbf{Mod}}}

\def\qQ{\mathbf{Q}}
\def\wA{\widetilde{A}}
\def\Hom{\mathrm{Hom}}
\def\id{\mathrm{id}}

\def\mc{\mathcal}

\def\cat#1{\mc{K}_2(A,#1)}  % homotopy category of duplexes

\def\drawing#1{\begin{center} \epsfig{file=#1} \end{center}}

\def\hsm{\hspace{0.03in}}

\def\doublemaprights#1#2#3#4{\raise3pt\hbox{$\mathop{\,\,\hbox to
      #1pt{\rightarrowfill}\kern-30pt\lower3.95pt\hbox to
      #2pt{\rightarrowfill}\,\,}\limits_{#3}^{#4}$}}

\newtheorem{theorem}{Theorem}
\newtheorem{prop}{Proposition}
\newtheorem{lem}{Lemma}
\newtheorem{corollary}{Corollary}

\numberwithin{equation}{section}
\numberwithin{lem}{section}
\numberwithin{prop}{section}

\newtheorem{definition}{Definition}

\newcommand{\oplusop}[1]{{\mathop{\oplus}\limits_{#1}}}

\newcommand{\bigoplusop}[1]{{\mathop{\bigoplus}\limits_{#1}}}

\begin{document}
\title{Homological realization of Nakajima varieties and 
 Weyl group actions}
\author{Igor Frenkel, Mikhail Khovanov and Olivier Schiffmann}
\maketitle
\baselineskip 14pt
\noindent

\tableofcontents

\vspace{0.2in}

\centerline{\textbf{Introduction}}

A geometric approach to representation theory of
Kac-Moody algebras was given by Nakajima in the groundbreaking work
\cite{Nak1}, which was a culmination of a series of remarkable 
discoveries discussed in the introduction to his paper. For any 
simply laced Kac-Moody algebra $\g$ with triples of generators 
$(e_a,h_a,f_a)_{a \in I}$ indexed by a finite
set $I$, Nakajima constructed a family of complex varieties
$\mathcal{M}_{\zeta}(\mathbf{v},\mathbf{w})$, where 
$\mathbf{v},\mathbf{w}
\in \N^I$ and $\zeta \in \RR^3 \otimes \RR^I$ such that, 
for any \textit{generic} $\zeta,$ 
$$\mathrm{dim\;}H^{mid}(\mathcal{M}_{\zeta}(\mathbf{v},\mathbf{w}))=
\mathrm{dim}\;L_{\lambda}[\lambda-\alpha],$$
where $L_{\lambda}$ is the integrable highest weight $\g$-module 
of highest weight $\lambda$, where $L_{\lambda}[\lambda-\alpha]$ is the 
corresponding weight space, and where
$$\alpha=\sum_a v_a\alpha_a,\qquad \lambda=\sum_a w_a\omega_a,$$
for $(\alpha_a)_{a\in I}, (\omega_a)_{a \in I}$ the set of simple roots 
and fundamental weights respectively. Nakajima realized 
the action of the generators $(e_a,h_a,f_a)_{a\in I}$ in a geometric way, 
the contravariant form on $L_{\lambda}$, and
indicated the geometric meaning of the Weyl group action on weight 
spaces
$$\sigma: L_{\lambda}[\lambda-\alpha] \stackrel{\sim}{\to} 
 L_{\lambda}[\sigma(\lambda-\alpha)].$$
The action of the Weyl group on the underlying quiver varieties was 
further developed by Lusztig \cite{LusWeyl}, Maffei \cite{Maffei}, 
and Nakajima \cite{Nak2}.

\paragraph{}Nakajima varieties are defined in terms of certain data
attached to the Dynkin diagram $\mathbf{Q}$ of the Kac-Moody algebra $\g.$ 
This data can be viewed as a generalization
of the Atiyah-Drinfeld-Hitchin-Manin description of the instanton 
moduli spaces. Correspondingly, various structures of representation 
theory of Kac-Moody algebra $\g$ including the action of the 
Weyl group on highest weight modules were described in terms of 
this linear data. We will review the original constructions 
relevant to our present work in Section~1. 

\vspace{.1in}

In this paper, we interpret the data and, consequently, Nakajima 
varieties, via differential graded modules over a 
finite-dimensional quotient $A(\qQ)$ of the double path 
algebra of $\qQ$. Defining relations in this quotient algebra depend 
on the choice of orientation $\epsilon$ of $\qQ,$ but different 
orientations produce isomorphic algebras. The algebra $A(\qQ)$ 
is the quadratic dual of the preprojective algebra of the 
(oriented) graph $\qQ$ and has a Frobenius structure. 
If $\qQ$ is bipartite (see Section 5), $A(\qQ)$ is isomorphic 
to the zigzag algebra of $\qQ$ studied in \cite{HK}. Let us 
denote $A(\qQ)$ simply by $A.$ 

Our realization of the Nakajima varieties allows us to view 
their theory in the context of homological algebra. In 
particular, simple $A$-modules $S_a$ and projective 
$A$-modules $P_a$ are the basic building blocks in our picture. 
Furthermore, we give a natural 
interpretation of the Weyl group action on Nakajima varieties, with 
the simple reflection $s_a$ acting as homological "addition" and 
 "substraction" of the projective module $P_a.$ This Weyl group 
action comes from a modification of the braid group action in 
the derived category of $A$-modules, see  \cite{KS}, \cite{ST}, 
\cite{RZ}, \cite{HK}.

\paragraph{} The theory of Nakajima varieties also suggest 
certain generalizations of some classical notions and results 
in homological algebra. In fact, our constructions depend in 
an essential way on the value of the 
parameter $\zeta=(\zeta_\RR,\zeta_\C)$, where $\zeta_\RR \in \RR^I$ and 
$\zeta_\C \in \C^I$.
If $\zeta_\C=0$ then we use the standard theory of differential graded 
algebras and complexes of modules over them. On the other hand, 
when $\zeta_\C \neq 0$  we are led to consider modules over $A$ equipped
with a generalized differential $d,$  which is a degree $1$ map 
satisfying
$$d^2=c$$ 
for a suitable central element $c \in A$. We introduce categories 
of $(A,c)$-complexes in Section 2. 

To define the Weyl group action we are then forced  into 2-periodic 
generalized complexes, which we call \textit{duplexes.} The theory of 
duplexes, which we outline in Section 3, can be developed in 
parallel with some classical results of homological 
algebra and seems to be worthy of a deeper, independent study.  

In Section~4 we consider another way of deforming the derived category of 
$A$-modules : by a theorem of Happel \cite{Happel}, this derived
category is equivalent to the stable category of graded modules over the algebra
$A \sotimes \C[d]/d^2$, where $\sotimes$ denotes the super tensor product. 
We introduce the stable category of
$A \sotimes \C[d]/\langle d^2-c \rangle$ and relate this to the category
of $(A,c)$-duplexes.

\paragraph{}Our interpretation of Nakajima varieties 
$\mathcal{M}_\zeta(\mathbf{v},\mathbf{w})$ is given in Section 5, 
and is set-theoretic, consisting of a bijection between the set 
of points of $\mathcal{M}_\zeta(\mathbf{v},\mathbf{w})$ and certain 
isomorphism classes of differential $A$-modules.
Let $(V_a)_{a \in I}$, $(W_a)_{a\in I}$ denote 
collections of vector spaces of dimensions $(v_a)_{a\in I}$ and $(w_a)_{a\in I}$ 
respectively.
We consider $\Z/2\Z$-graded $A$-modules $M$ equipped with a 
generalized differential $d$ such that $d^2=c$, where $c$ is 
$\zeta_\C$, viewed as a central element of $A,$ and 
\begin{equation}\label{E:intro1}
M \cong \bigoplus_a (P_a \otimes V_a \oplus S_a [-1]\otimes W_a),
\end{equation}
where $[-1]$ denotes a grading shift, and the isomorphism 
is that of $A$-modules. When $\zeta_\C$ is not generic 
(for example $\zeta_\C=0$), we add an irreducibility condition with
respect to $d$, analogous to a stability condition. We show that 
the set of classes of pairs $(M,d)$ as above is in a natural 
bijection with the set of $G_W$-orbits of the quiver variety 
$\mathcal{M}_\zeta(\mathbf{v},\mathbf{w})$, where $G_W=\prod_a GL(W_a)$. The
full variety $\mathcal{M}_\zeta(\mathbf{v},\mathbf{w})$ is obtained by fixing in addition
an isomorphism (framing) $W_a \simeq \mathrm{Hom}_A(S_a[-1], M)$ for $a \in I$.

\paragraph{}To obtain a realization of the Weyl group action, we 
consider a duplex of $A$-bimodules  $C_{a,x}$
associated to a vertex $a$ of $Q$ :
$$ \to P_a \otimes \hsm _a P \to A\to$$
where $x\in \C$ and one of the bimodule maps depends on $x.$ 
We show in Section 6 that duplexes $C_{a,x}$ are invertible in the 
homotopy category (as well as in the stable category), 
$$
C_{a,x} \otimes_A C_{a,-x} \cong A,
$$
and satisfy Yang-Baxter relations:
$$
C_{a,x} \otimes C_{b,y} \cong C_{b,y} \otimes C_{a,x}
$$
for any two vertices $a$ and $b$ which are not joined by an edge, and 
$$ C_{a,x} \otimes C_{b,x+y} \otimes C_{a,y}
\cong C_{b,y} \otimes C_{a,x + y} \otimes C_{b,x}$$
for $a$ and $b$ joined by a single edge.
In the limit $x\to 0$ our duplexes degenerate into the ones used to 
categorify the Burau representation of the braid group, see 
\cite{KS}. 

Points of Nakajima varieties 
$\mathcal{M}_{\zeta}(\mathbf{v},\mathbf{w})$
can be identified with certain isomorphism classes of $A$-duplexes.
The functor of the tensor product with the bimodule duplex $C_{a,x}$ 
acts in categories of $A$-duplexes and restricts to a bijection 
$$\mathcal{R}_a: \mathcal{M}_{\zeta}(\mathbf{v},\mathbf{w}) 
 \stackrel{\sim}{\to} \
\mathcal{M}_{s_a\cdot\zeta}(s_a\cdot (\mathbf{v},\mathbf{w})),$$
where $s_a$ is a simple reflection. We show in Section 7 that our 
reflection maps coincide with those in \cite{LusWeyl}, \cite{Maffei}, 
and \cite{Nak2}.

\paragraph{} The class of Kac-Moody algebras associated to affine 
Dynkin diagrams plays a very special role in representation theory. 
The corresponding class of quivers is also distinguished in the 
Nakajima theory since it appears in the study of the instanton 
moduli spaces. It is well-known that simply-laced affine Dynkin 
diagrams are in a bijection with finite subgroups $\Gamma\subset 
SL(2,\C)$ and it is natural to reformulate Nakajima's work 
entirely in terms of these finite groups. In Section 8 we recast 
our realization of Nakajima varieties in this light. We replace 
$A$ by the Morita equivalent algebra 
$A_\Gamma = \Lambda \rho \otimes \C[\Gamma],$ where $\rho$ 
is the natural 2-dimensional representation of $\Gamma.$ To a 
collection of vector spaces $(V_a)_{a\in I}$ and $(W_a)_{a\in I}$ we now 
associate two $\Gamma$-modules
 $$\mathbb{V}=\bigoplus_a \rho_a \otimes V_a, \qquad 
 \mathbb{W}=\bigoplus_a \rho_a \otimes W_a,$$
where $(\rho_a)_{a\in I}$ is the set of all irreducible representations 
of $\Gamma$. The module $M$ in (\ref{E:intro1}) is replaced 
by the following 
 \begin{equation}\label{U:intro1} 
  M = \Lambda\rho \otimes \mathbb{V} \oplus \mathbb{W}[-1].
 \end{equation}  
The latter can be viewed as a module over 
$\widetilde{A}_{\Gamma,c}=A_\Gamma \sotimes \C[d]/\langle d^2-
 c\rangle$, where $c$ is a central element of $A_\Gamma$ 
which depends linearly on $\zeta_{\C}$.  It turns out that 
we can characterize modules of the form (\ref{U:intro1}) 
as a certain class of elements of the stable category of $\Z/2\Z$-graded 
$\widetilde{A}_{\Gamma,c}$-modules 
$\underline{\Mod}_2(\widetilde{A}_{\Gamma,c}).$
This yields a realization of Nakajima varieties 
via isomorphism classes of pairs $(M,u),$ where $M$ is an 
object of $\underline{\Mod}_2(\widetilde{A}_{\Gamma,c})$
and $u:\mathbb{W}\to R(M)$ is a fixed isomorphism, with $R(M)$ 
being the restriction of the module $M$ to $A_{\Gamma}.$ 

\vspace{.1in}

The relation between this realization of the Nakajima varieties by means of the stable 
category $\underline{\Mod}_2(\widetilde{A}_{\Gamma,c})$ and
the realization as a moduli space of $\Gamma$-equivariant 
torsion-free sheaves on a  \textit{noncommutative} $\mathbb{P}^2$ with 
fixed framing at infinity (presented in \cite{BGK}) is explained in Section~9. 
In fact, one may view such torsion-free sheaves as modules over the  
algebra Koszul dual to $\widetilde{A}_{\Gamma,c},$ see 
\cite[Appendix B]{BGK}.
Thus our construction illustrates a noncommutative version 
of the classical theorem of \cite{BGG} claiming
that when $\Gamma=\{e\}$ and $c=0$ (the case of commutative 
$\mathbb{P}^2$), 
the derived category of coherent sheaves over $\mathbb{P}^2$
is equivalent to the stable category of $\widetilde{A}_{\Gamma,c}$. 
Finally,
we reformulate the action of the Weyl group on quiver varieties via
 certain natural duplexes of $A_\Gamma$-bimodules.

\paragraph{}We believe that the realization of the Nakajima 
varieties and Weyl group actions by means of natural categorical 
constructions presented in this paper is only a first step in a more 
general program of recasting the
Nakajima geometric approach to representation theory of Kac-Moody 
algebras in
terms of canonical structures of homological algebra. We hope that 
the emerging interaction between the two areas will be beneficial to 
both subjects. Below we will make a few remarks about further 
developments of both areas inspired by our constructions. 

\paragraph{}As we mentioned in the beginning of Introduction, the 
Nakajima varieties encode the structure of the integrable highest 
weight modules $L_{\lambda}.$ The latter modules possess rich 
structures associated with the corresponding Weyl group $\mathbf{W}$. On 
one hand, each module $L_{\lambda}$ contains a family of 
Demazure submodules $L_{\lambda,w}$ defined for any $w\in \mathbf{W}$. On 
the other, each module $L_{\lambda}$ admits the BGG resolution 
by Verma modules $V_{w\cdot \lambda},$ where again $w\in \mathbf{W}$ 
(see e.g. \cite{Ku} for a review of both constructions). It is 
natural to expect that our realization of the action of $\mathbf{W}$ on 
Nakajima varieties should lead to a transparent geometric 
construction of the Demazure modules as well as the BGG resolution. 

Concerning the applications to homological algebra, the 
example studied in this paper already suggests the following 
generalizations of algebraic structures :

\begin{align*}
\text{Abelian} \qquad &\longrightarrow \qquad \text{Triangulated}\\
\Z\text{-graded}\qquad &\longrightarrow \qquad\Z/2\Z\text{-graded}\\
d^2=0 \qquad &\longrightarrow \qquad d^2=c
\end{align*}

Starting with 
an abelian category of modules over a ring, we pass to the 
triangulated category of complexes (the top arrow in the diagram). 
This arrow is the familiar advancement from the classical 
theory of modules over a ring to homological algebra. The 
bottom arrows refer to two more recent developments where 

\begin{itemize} 
\item  one gains from working with periodic triangulated 
categories, those with $[2k]\cong\mathrm{Id}$ for some $k$
(case $k=1$ seems especially important),  
\item differential modules acquire curvature (the square of 
 differential is no longer zero). 
\end{itemize} 
 
Both transformations are natural from the deformation theory 
viewpoint. When the cohomology ring of a symplectic manifold 
is deformed to the quantum cohomology ring, its $\Z$-grading 
collapses to a $\Z/{2k}\Z$-grading, where $k$ is the minimal 
Chern number (see [MS, Section 1.7]),  
and the shift functor in the $A_{\infty}$ triangulated 
Fukaya-Floer category of the manifold is periodic, 
$[2k]\cong\mathrm{Id}.$ Deforming $d^2=0$ to $d^2=c$ is 
as legitimate as deforming the ring structure, when 
one is describing all deformations of the homotopy category 
of modules over  a (graded) ring. We should also mention 
the paper of L.Peng and J.Xiao \cite{PX}, where 2-periodic 
triangulated categories appear in relation to Hall algebras.  

\vspace{0.1in}

{\bf Acknowledgments:} M.K. would like to thank Igor Burban for interesting 
discussions about periodic triangulated categories. While working 
on this paper, I.F. and M.K. were partially supported by NSF 
grants DMS-0070551 and DMS-0104139. 

\section{Nakajima varieties}

\begin{figure} [htb] \drawing{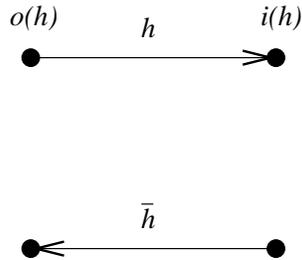}\caption{Two orientations of 
 an edge} \label{edge} 
 \end{figure}

\paragraph{} We recall the definition of Nakajima quiver varieties.
Let $\qQ=(I, E)$ be an arbitrary finite graph with $I$ the set of 
vertices and $E$ the set of edges. We allow $\qQ$ to have loops 
and multiple edges. Let $H$ be the set of oriented edges of this graph 
(thus $H$ is "twice as large" as $E$). For any $h \in H$ we denote by 
$o(h)$ and $i(h)$ the outgoing and incoming vertices of $h,$ 
respectively, and by $\overline{h}$ the edge $h$ with the opposite 
orientation, see figure~\ref{edge}.

Let $(\mathbf{v},\mathbf{w})\in \N^I \times \N^I$ with
$\mathbf{v}=(v_a)_{a \in I}$ and
$\mathbf{w}=(w_a)_{a \in I}$. Fix some $I$-graded $\C$-vector spaces
$V=\bigoplus V_a$ and $W=\bigoplus W_a$ such that
$\mathrm{dim}\;V=\mathbf{v}$ and $\mathrm{dim}\;W=\mathbf{w}$. Set
$$E(V,V)=\bigoplus_{h \in H} \mathrm{Hom}(V_{o(h)}, V_{i(h)}),$$
$$L(V,W)= \bigoplus_{a} \mathrm{Hom}(V_a, W_a),\qquad
L(W,V)=\bigoplus_{a} \mathrm{Hom}(W_a, V_a)$$
and
$$\mathbf{M}(\mathbf{v},\mathbf{w})=E(V,V) 
\oplus L(W,V) \oplus L(V,W).$$
An element of $\mathbf{M}(\mathbf{v}, \mathbf{w})$ will usually be 
denoted by its components $(B,i,j)$.

Let $\epsilon : H \to \{1,-1\}$ be any function satisfying
$\epsilon(h)+\epsilon(\overline{h})=0$ for all $h\in H.$ Such 
functions are 
in a bijection with orientations of $\qQ,$ the $\epsilon$-orientation 
consists of all edges $h$ with $\epsilon(h)=1.$ Consider the maps 
\begin{eqnarray*}
\mu_{\mathbb{C}} : 
 \mathbf{M}(\mathbf{v},\mathbf{w})&\lra & \bigoplus_a 
\mathfrak{gl}(V_a),\\
(B,i,j) &\lra  & 
\left( \sum_{o(h)=a} \epsilon(h)B_{\overline{h}}B_h+i_aj_a \right)_a,
\end{eqnarray*}
and 
\begin{eqnarray*}
\mu_\RR:\; \mathbf{M}(\mathbf{v},\mathbf{w}) &\lra &  \bigoplus_a
\mathfrak{u}(V_a),  \\
(B,i,j) &\mapsto &  
\frac{\sqrt{-1}}{2}\left( \sum_{o(h)=a} B_{\overline{h}}
B_{\overline{h}}^\star -B_{h}^\star  
B_{h} + i_a i_a^\star - j_a^\star j_a \right)_a.  
\end{eqnarray*}
In the above, $f^\star$ denotes the Hermitian adjoint of $f$. 
Following Nakajima, to $\zeta_{\R} \in \mathbb{R}^I$ we associate a central element 
$\boldsymbol{\zeta}_{\R}=\sum_a\frac{\sqrt{-1}\zeta_{\R,a}}{2}Id \in \bigoplus_a \mathfrak{u}(V_a)$,
and to $\zeta_{\C} \in \C^I$ we associate a central element
$\boldsymbol{\zeta}_{\C}=\bigoplus_a \zeta_{\C,a}Id \in \bigoplus_a \mathfrak{gl}(V_a)$.
 The group ${U}_V=\prod_a U(V_a)$ acts on 
$\mathbf{M}(\mathbf{v},\mathbf{w})$ by conjugation.
Finally, we put
$$\mathcal{M}_{\zeta}(\mathbf{v},\mathbf{w})
=(\mu_\RR \times \mu_\C)^{-1}(\boldsymbol{\zeta}_{\R},\boldsymbol{\zeta}_{\C})/{U}_V,$$
where $\zeta=(\zeta_{\R},\zeta_{\C})$.
Different choices of $\epsilon$ yield isomorphic varieties.

\paragraph{1.3.} When $\zeta_{\R} \in \Z^I$ 
there is also a purely complex-geometric
description of 
$\mathcal{M}_{\zeta}
(\mathbf{v},\mathbf{w})$. Note that
the group $G_V=\prod_a GL(V_a)$ acts on $\mathbf{M}(\mathbf{v},
\mathbf{w})$ by conjugation. To $\zeta_{\R}$ we
associate the character
$\chi_{\zeta_{\R}}:
\;G_V \to \C^*,\;(g_a)_a\mapsto\prod_a det(g_a)^{\zeta_{\R,_a}}$. 
Following Nakajima \cite{Nak1}, we have
$$\mathcal{M}_{\zeta}
(\mathbf{v},\mathbf{w})\simeq \mathrm{Proj}\;
\big(\bigoplus_{n \geq 0} A^n_{\zeta_{\R}}(\mathbf{v},\mathbf{w})\big),$$
where
\begin{equation*}
\begin{split}
A_{\zeta_{\R}}^n&(\mathbf{v},\mathbf{w})\\
&=\{f \in \C[\mu_{\C}^{-1}(\boldsymbol{\zeta}_{\C})]
\;|\;f(g\cdot (B,i,j))=
\chi_{\zeta_{\R}}(g)^nf((B,i,j))\quad \forall\;g \in G_V\}.
\end{split}
\end{equation*}

\paragraph{}There is an open subset
$\mathbf{M}^{ss}_{\zeta}(\mathbf{v},\mathbf{w})
\subset \mu_{\C}^{-1}(\boldsymbol{\zeta}_{\C})$ of \textit{semistable points} such that
$\mathcal{M}_{\zeta}(\mathbf{v},\mathbf{w})=
\mathbf{M}^{ss}_{\zeta}(\mathbf{v},\mathbf{w})
//G_V$ (categorical quotient). We will mainly be interested in the following
two cases :
\begin{enumerate}
\item[i)] $\zeta_{\R} \in (\mathbb{N}^+)^I$ and $\zeta_{\C}=0$. In this case
$(B,i,j) \in \mathbf{M}^{ss}_{\zeta_{\C}}(\mathbf{v},
\mathbf{w})$ if the following condition is satisfied: the only (graded)
$B$-invariant subspace of $V$ contained in $\mathrm{Ker}\;{j}$ is $\{0\}$.
\item[ii)] $\zeta_{\R}$ is arbitrary and $\zeta_{\C}$ satisfies 
the
following genericity condition: for every $n_1,\ldots,n_k \in \Z$, we have
$\sum_a n_a \zeta_{\C,a} = 0 \Rightarrow n_a=0 $ for all $a$. In this
case, all points $(B,i,j)$ in $\mu_{\C}^{-1}(\boldsymbol{\zeta}_{\C})$ are semistable (and in fact
all points $(B,i,j)$ in $\mu_{\C}^{-1}(\boldsymbol{\zeta}_{\C})$ automatically satisfy the condition
in i)). 
\end{enumerate}
In the two above cases, $G_V$ acts freely on
$\mathbf{M}^{ss}_{\zeta}(\mathbf{v},\mathbf{w})$
(see e.g \cite{Nak1}), so that the categorical quotients are actually geometric
(smooth) quotients. Note also that, in the case ii), the variety
$\mathcal{M}_{\zeta}(\mathbf{v},\mathbf{w})$
is actually independent of $\zeta_{\R}$.

\paragraph{1.4.} To the graph $(I,E)$ we associate a symmetric $|I| \times |I|$ Borcherds matrix
$A=(a_{ij})$ with 
$$a_{ij}=2\delta_{ij}-\#\{h \in H\;| i(h)=i, o(h)=j\}.$$
Let $I^{re} \subset I$ be the set of all loopless vertices (characterized by the relation
$a_{ii}=2$).
To $A$ corresponds a Borcherds algebra (or Generalized Kac-Moody algebra) $\g$
(see \cite{Bo}). Let us fix a Cartan decomposition $\g=\n_- \oplus \h \oplus
\n_+$ and let $\alpha_a$ and $\omega_a$ (resp. $\alpha^\vee$ and 
$\omega_a^\vee$) 
stand for the simple root and
fundamental weight (resp. simple coroot and fundamental coweight)
associated to a vertex $a$. We put 
$$Q=\bigoplus_s \Z \alpha_a, \qquad
 P=\bigoplus_a \Z \omega_a.$$
We also set $Q^{\vee} =\bigoplus_a \Z \alpha_a^{\vee}$ and 
$P^\vee=\bigoplus_a \Z \omega_a^\vee$
and we denote by $\langle ,\rangle$
the natural pairing between $\h$ and $\h^*$.
We consider $\mathbf{v},
\mathbf{w},\zeta_{\R}$ and $\zeta_{\C}$ as elements of $\h^*$
via the identifications
\begin{equation*}
\mathbf{v} \mapsto \sum_a v_a \alpha_a,\qquad
\mathbf{w} \mapsto \sum_a w_a \omega_a,\qquad
\zeta_{\R} \mapsto
\sum_a \zeta_{\R,a}\omega_a,\qquad \zeta_{\C} \mapsto
\sum_a \zeta_{\C,a} \omega_a.
\end{equation*}

\paragraph{1.5.} We say that the parameter $\zeta$ is \textit{generic} when
\begin{equation}\label{E:2}
\mathrm{For\;every\;}\nu \in P^\vee,\;\mathrm{we\;have\;}
<\nu,\zeta_{\R}>\neq 0
\;\mathrm{or}\;<\nu,\zeta_{\C}>\neq 0
\end{equation}

The variety $\mathcal{M}_{\zeta}
(\mathbf{v},\mathbf{w})$ is smooth whenever
$\zeta$ is generic. Moreover (for fixed
$\mathbf{v}$ and $\mathbf{w}$) the varieties corresponding to generic
parameters are all diffeomorphic
(\cite{Nak1}, Corollary 4.2).

\paragraph{1.6.} The Weyl group $\mathbf{W}$ of $\g$
is defined to be the subgroup of $\mathrm{Aut}(\h^*)$ generated by reflections
$$s_a: \alpha \mapsto \alpha-\langle \alpha, \alpha_a^\vee \rangle \alpha_a$$
for $a \in I^{re}$. Note that $\mathbf{W}$ acts on $Q$. The dual action
on $Q^\vee$ is given by
\begin{equation}\label{E:Waction}
s_a: \alpha^\vee \mapsto \alpha^\vee - \langle \alpha_a, \alpha^\vee \rangle \alpha^\vee_a.
\end{equation}
Moreover, $\langle\,,\, \rangle$ induces a perfect pairing between $P$ and $Q^\vee$, and thus (\ref{E:Waction})
gives rise to an action of $\mathbf{W}$ on $P$ by duality.

With this convention, $\mathbf{W}$ acts on $\zeta_{\R}$ and
$\zeta_{\C}$ via the identifications in 1.4. We also define an
action on pairs $(\mathbf{v},\mathbf{w})$ by $\sigma \cdot (\mathbf{v},
\mathbf{w}):= (\sigma(\mathbf{v}-\mathbf{w}) + \mathbf{w},\mathbf{w})$.
Following \cite{LusWeyl} (see also \cite{Nak1}), Maffei \cite{Maffei}
defined, for \textit{generic} $\zeta$,
isomorphisms
$$\kappa_\sigma:\;\mathcal{M}_{\zeta}(\mathbf{v},
\mathbf{w}) \overset{\sim}{\to}
\mathcal{M}_{\sigma\cdot\zeta}
(\sigma\cdot(\mathbf{v},\mathbf{w})).$$

\paragraph{}Let us consider the situation when
$\zeta_{\C}$ is generic as in Section 1.3. ii) (so that
$\zeta$ is generic). In that situation, it is more convenient to
use the purely complex description of the quiver variety given in Section~1.3. In the case of a 
simple reflection
$s_a$, the construction is as  follows. Let $(B,i,j) \in
\mathbf{M}^{ss}_{\zeta} (\mathbf{v},\mathbf{w})$.
Define vector spaces $V'_k, W'_k$ by $V'_k=V_k$ if $k\neq a$,
$$V'_a=(W_a \oplus \bigoplus_{o(h)=a}V_{i(h)})/(j_a +\sum_{o(h)=a}B_h)V_a$$
and
$W'_k=W_k$ for all $k$. Let $Z$ be the set of all $(B',i',j') \in
\mathbf{M}_{s_a \zeta}
(\mathbf{v}',\mathbf{w}')$ such that :
\begin{enumerate}
\item[i)] $B'_h=B_h$ if $i(h) \neq a$ and $o(h) \neq a$,
\item[ii)]$i_k=i'_k$ and $j'_k=j_k$ if $k \neq a$,
\item[iii)] Set
$$x_a=j_a \oplus \bigoplus_{o(h)=a}B_{h} : V_a \to W_a \oplus \bigoplus_{o(h)=a}
V_{i(h)},$$
$$y_a=i_a \oplus \bigoplus_{i(h)=a}\epsilon(\overline{h})B_{h} : W_a \oplus 
\bigoplus_{i(h)=a}
V_{o(h)}\to V_a$$
and define $x'_a$ and $y'_a$ in a similar fashion. The sequence
$$0 \to V_a \stackrel{x_a}{\lra} W_a \oplus \bigoplus_{o(h)=a}V_{i(h)} 
\stackrel{y'_b}{\lra}
V'_a \to 0$$
is exact, and $x_ay_a=x'_ay'_a-\lambda_a Id$.
\end{enumerate}
Then (see \cite{Maffei}), $Z$ is a principal $GL(V_a)$-homogeneous space. Thus it corresponds
to a unique point $\kappa_{s_a}(B,i,j) \in \mathcal{M}_{s_a\zeta}(\mathbf{v}',\mathbf{w}')=
\mathcal{M}_{s_a\zeta} (s_a\cdot(\mathbf{v},\mathbf{w}))$.

\section{Categories of $(A,c)$-complexes}

\paragraph{2.1.} Let $A$ be a $\Z$-graded ring and $c$ a central element of
$A$ of
degree $2.$ We denote by $Z_2(A)$ the degree $2$ summand of the center of
$A,$ so that $c\in Z_2(A).$

\begin{definition} A left $(A,c)$-complex is a $\Z$-graded left
$A$-module $M$ together with a degree one map $d: M \to M$ such that
  \begin{equation*}
    d^2=c, 
  \end{equation*}
and $d$ (super)commutes with the action of $A$:
   \begin{equation*}
     d (a m ) = (-1)^{|a|} a d(m), \hspace{0.2in} a\in A, m \in M.
   \end{equation*}
\end{definition}

A morphism of $(A,c)$-complexes is a degree zero morphism of $A$-modules
which commutes with $d$. The category $\mathbf{Com}(A,c)$
of left $(A,c)$-complexes is abelian. The translation functor $[1]$ in the
category $\mathbf{Com}(A,c)$
is defined by
   \begin{equation*}
      (M[1])^i = M^{i+1}, \hspace{0.1in} d_{[1]} = - d,
   \end{equation*}
  and the $A$-module structure on $M[1]$ is
  \begin{equation*}
   a\circ m = (-1)^{|a|} a m .
  \end{equation*}

Alternatively, we may define a $\Z$-graded algebra
$\widetilde{A}_c=A \sotimes \C[d]/(d^2-c)$ with $da=(-1)^{|a|}ad$ and
$deg(d)=1$. Then $\mathbf{Com}(A,c)$ is nothing but the category
$\mathbf{Mod}(\widetilde{A}_c)$ of graded left $\widetilde{A}_c$-modules.

\paragraph{2.2.} Let $M$ and $N$ be left $(A,c)$-complexes.
Given a morphism of graded $A$-modules $h: M \to N[-1],$ the map
$f = h d_M + d_N h $
is a morphism $M\to N$ of $(A,c)$-complexes. Any such morphism is
called \emph{null-homotopic}. The following result is clear.

\begin{prop} Null-homotopic morphisms form a 2-sided ideal in
the category $\mathbf{Com}(A,c)$.
\end{prop}

We say that morphisms $f,g: M\to N$  are homotopic and write $f \sim g$
if $f-g$ is null-homotopic. Define the homotopy category 
$\mc{K}(A,c)$ as follows. Objects are 
$(A,c)$-complexes and for any two $(A,c)$-complexes $M$ and $N$ we put
$$\mbox{Hom}_{\mc{K}(A,c)}(M,N)=
\mbox{Hom}_{\mathbf{Com}(A,c)}(M,N)/ \sim .$$

Categories $\mathbf{Com}(A,c),$ as well as $\mc{K}(A,c),$ for various  
$c\in Z_2(A),$ might have common objects. If $M\in \mathbf{Com}(A,c)$ and 
$c' M =0$ for some $c'\in Z_2(A)$ then $M\in \mathbf{Com}(A,c+c').$

\paragraph{2.3.}\emph{Tensor product of left and right $(A,c)$-complexes.}

\vspace{0.1in}

If $M$ is a right graded $A$-module and $N$ a left graded $A$-module,
the tensor product $M\otimes_A N$ is a graded abelian group.
If $M$ is a right $(A,c)$-complex and $N$ a left $(A,-c)$-complex,
then $M\otimes_A N$ is a complex of graded abelian groups
with the differential
   \begin{equation*}
    d(m \otimes n) = d m \otimes n + (-1)^{|m|} m \otimes dn,
   \end{equation*}
since
   \begin{equation*}
     d^2(m \otimes n) = d^2 m \otimes n + m \otimes d^2 n = m c
     \otimes n + m \otimes (-c) n = 0.
   \end{equation*}

\paragraph{2.4.}\emph{Bimodules}

\vspace{0.1in}

Let $c_0,c_1\in Z_2(A).$ An $(A,c_0,c_1)$-complex is a graded $A$-bimodule
$N$ together with a degree one map $d: N\to N$ such that
  $d^2= l_{c_0}  + r_{c_1},$
(where $l_{c_0}$ is left multiplication by $c_0$ and $r_{c_1}$ is right 
multiplication by $c_1$), $d$ (super)commutes with the left action of $A$:
  \begin{equation*}
     d (a n ) = (-1)^{|a|} a \hsm d n, \hspace{0.2in} a\in A, n \in N,
   \end{equation*}
  and commutes with the right action of $A.$ 

If $M$ is a left $(A,-c_1)$-complex, the tensor product $N\otimes_A M$ is
a left $(A,c_0)$-complex. Thus, the tensor product with $N$ is a functor
from $\mathbf{Com}(A,-c_1)$ to $\mathbf{Com}(A,c_0),$ 
and from $\mc{K}(A,-c_1)$ to $\mc{K}(A,c_0).$

\section{Categories of duplexes}

\paragraph{3.1.}Let $A$ be a $\Z/2\Z$-graded ring, $A= A_0\oplus A_1,$ and 
$c$ a degree zero central element of $A.$ 
A \emph{duplex} over $(A,c)$ is a $\Z/2\Z$-graded $A$-module 
$M=M^0\oplus M^1$ with a generalized differential 
\begin{equation} \label{diff-dupl}
     M^0\stackrel{d}{\lra} M^1 \stackrel{d}{\lra} M^0
\end{equation}
 which supercommutes with the action of $A$ and satisfies $d^2(m)=cm$ for 
 all $m\in M.$ A duplex over $(A,0)$ is simply a 2-periodic complex of 
 $A$-modules. An $(A,c)$-duplex for $c\not= 0$ can be viewed as a
 2-periodic "complex" with the differential satisfying $d^2=c$ rather
 than $d^2=0.$

A homomorphism $f:M\to N$ of $(A,c)$-duplexes is a degree zero $A$-module 
map that commutes with the differentials:   

\begin{equation*}
    \begin{CD}
      M^0 @>{d}>> M^1 @>{d}>> M^0   \\
      @V{f^0}VV      @V{f^1}VV   @V{f^0}VV \\
      N^0 @>{d}>> N^1 @>{d}>> N^0
    \end{CD}
\end{equation*}
We denote by $\mathbf{Com}_2(A,c)$ the category of $(A,c)$-duplexes. 
The shift functor $[1]$ in this category is 2-periodic, $[2]\cong\mbox{Id}.$
The category $\mathbf{Com}_2(A,c)$ is abelian. We will often write duplexes 
in the form
\begin{equation*}
   \stackrel{d}{\lra} M^0\stackrel{d}{\lra} M^1 \stackrel{d}{\lra}
\end{equation*}

Given maps $f: M \to N$ and $g: N\to M$ of $(A,c)$-duplexes such 
that $fg= c_1 $ and $gf= c_1$ for some degree zero central element 
$c_1$ of $A,$ 
\begin{equation*}
    \begin{CD}
      M^0 @>{d_M}>> M^1 @>{d_M}>> M^0   \\
      @V{f^0}VV      @V{f^1}VV   @V{f^0}VV \\
      N^0 @>{d_N}>> N^1 @>{d_N}>> N^0   \\
      @V{g^0}VV      @V{g^1}VV   @V{g^0}VV \\
      M^0 @>{d_M}>> M^1 @>{d_M}>> M^0   \\
    \end{CD}
\end{equation*}
the \emph{cone} of $(f,g)$ is defined as the "total" 
 $(A,c+c_1)$-duplex of the above diagram, 
 \begin{equation*} 
  \stackrel{d}{\lra} M^0 \oplus N^1 \stackrel{d}{\lra} M^1 \oplus N^0 
  \stackrel{d}{\lra}
 \end{equation*} 
 with $d= d_M - d_N + f + g.$  

We will also use a less precise notation 
$\stackrel{g}{\lra} M \stackrel{f}{\lra} N \stackrel{g}{\lra} $ for 
the cone of $(f,g).$ 

\paragraph{}Note that we may once again think of $\mathbf{Com}_2(A,c)$
as the category of $\Z/2\Z$-graded left $\widetilde{A}_c$-modules, where
$\widetilde{A}_c$ is defined as in Section~2.

\paragraph{3.2.}
Given an $A$-homomorphism $h: M \to N[-1],$ the map $f = h d_M + d_N h$
is a morphism $M\to N$ of duplexes. We will say that morphisms $f,g:
M\to N$ are homotopic if $f- g = h d_M + d_N h $ for some $h.$ The following
is straightforward.

\begin{prop} Null-homotopic morphisms form a 2-sided
ideal in the category $\mathbf{Com}_2(A,c).$
\end{prop}

We call the quotient category of $\mathbf{Com}_2(A,c)$ by this ideal the
\emph{homotopy category of $(A,c)$-duplexes} and denote it by
$\cat{c}.$

\vspace{.1in}

\noindent
\emph{Example:} For an $A$-module $M$ let $M_{c,1}$ be the duplex
  \begin{equation} \label{short}
   \stackrel{1}{\lra} M \stackrel{c}{\lra} M \stackrel{1}{\lra}, 
  \end{equation}
which is the cone of $(c,1).$ 
The identity morphism of $M_{c,1}$  is null-homotopic, and $M_{c,1}$
is isomorphic to the zero object in the homotopy category of duplexes.

\vspace{0.1in}

\noindent
\textbf{Remark:} If $c$ is invertible, any $(A,c)$-duplex
is trivial in the homotopy category, and the category $\cat{c}$ is trivial.
The case of non-invertible $c$ is more interesting. If $A$ is
artinian, any element of $A$ is either invertible or nilpotent, and
the only nontrivial case is that of nilpotent $c.$

\vspace{0.1in}

\begin{prop} \label{quot-by-injective} If $M$ is an $(A,c)$-duplex 
and $I$ an injective $\Z/2\Z$-graded $A$-submodule of $M$ such 
that $d$ is injective on $I$ and $I\cap dI=0,$ then $M$ is 
isomorphic in the homotopy category to its quotient by the subduplex 
generated by $I$: 
\begin{equation*}
M \cong \{ \lra M^0/(I^0\oplus d(I^1)) \lra M^1/(I^1\oplus d(I^0)) \lra \}. 
\end{equation*}
\end{prop}

\noindent
The proof is again straightforward. \qed

\paragraph{3.3.} \emph{Duplexes of bimodules}\\
Tensor product with a duplex $N$ of
$A$-bimodules such that $d^2 = l_{c_0} + r_{c_1}$
is a functor from $\mathbf{Com}_2(A,-c_1)$ to $\mathbf{Com}_2(A,c_0)$ and from
$\mc{K}_2(A,-c_1)$ to $\mc{K}_2(A,c_0).$

\section{Stable categories}

\paragraph{4.1.} Let $A$ be a $\Z$-graded ring and $c \in Z_2(A)$ a degree two central
element.
Let $\UMod (A) $ denote the stable category of
$\Z$-graded $A$-modules (see e.g \cite{Happel}).
Its objects are $\Z$-graded $A$-modules
and for any modules $M$ and $N$ we have
$$\mathrm{Hom}_{\UMod (A)}(M,N)=\mathrm{Hom}_{A}(M,N)/I$$
where $I$ is the ideal of all morphims $f: M \to N$ which admit a factorisation
$f=g\circ h$ where $h: M \to P, \;g:P\to N$ for some projective module $P$.
In particular, an object $M$ of $\UMod (A)$ is isomorphic to
the zero object if and only if it is projective as an $A$-module.
We define the stable category $\UMod(\wA_c)$ in a same way. Since
$\wA_c$ is projective (in fact free) as an $A$-module, there is a
natural restriction functor $R: \UMod(\wA_c) \to \UMod(A)$.

\paragraph{}There is a canonical functor $\Phi: \mathbf{Mod}(\wA_c)\simeq\mathbf{Com}(A,c)
\to \mathcal{K}(A,c)$.
\begin{lem} For any projective $\wA_c$-module $P$ we have $\Phi(P)=0$.
\end{lem}
\noindent
\textit{Proof.} Let $\{P'_i\}$ be the collection of indecomposable
projective $A$-modules. It is easy to check that $P_i=\wA_c \otimes_A P'_i$ is
an indecomposable projective $\wA_c$-module, and hence that $\{P_i\}$
forms the complete collection of indecomposable projectives for $\wA_c$. But
$P_i=P'_i \oplus P'_i[-1]$ as $A$-module and $d: P'_i \stackrel{\sim}{\to} P'_i[-1]$,
so that $P_i$ is homotopic to zero as an $(A,c)$-complex.\qed

\paragraph{}We deduce that the functor $\Phi$ admits a factorization
\begin{equation}\label{E:fonct}
\mathbf{Mod}(\wA_c)\simeq\mathbf{Com}(A,c) \stackrel{\Phi_1}{\longrightarrow}
\UMod(\wA_c) \stackrel{\Phi_2}{\longrightarrow} \mathcal{K}(A,c).
\end{equation}
Similar results hold for $\mathbf{Com}_2(A,c)$, $\UMod_2(\wA_c)$ and $\mathcal{K}_2(A,c)$
if $A$ is a $\Z/2\Z$-graded ring.

\paragraph{4.2.} Proposition~\ref{quot-by-injective} admits the following straightforward
generalization.
\begin{prop}\label{quot-by-projinj}
If $M$ is an $(A,c)$-duplex 
and $I$ an injective and projective $\Z/2\Z$-graded $A$-submodule of $M$ such 
that $d$ is injective on $I$ and $I\cap dI=0,$ then $M$ is 
isomorphic in the stable category $\UMod_2(\wA_c)$ to its quotient by the subduplex 
generated by $I$: 
\begin{equation*}
M \cong \{ \lra M^0/(I^0\oplus d(I^1)) \lra M^1/(I^1\oplus d(I^0)) \lra \}. 
\end{equation*}
\end{prop}

\section{Homological realization of Nakajima varieties}

\paragraph{5.1.} Let $\qQ=(I,E), H$ and $\epsilon$ be as in Section 1.1. 
We can view $(I,H)$ as the oriented double of the unoriented 
graph $\qQ$. Consider the path algebra of $(I,H).$ Note that in this 
algebra the product $h h'$ of two length one paths is nonzero if and only 
if $i(h)= o(h').$ 

Define the $\C$-algebra $A(\qQ)$ as the 
quotient of this path algebra by relations 
 \begin{enumerate}
  \item[i)] $ h h'=0$ if $h' \not= \overline{h},$   
  \item[ii)] $\epsilon(h) h \overline{h} =\epsilon(h') h'\overline{h'}$ 
  if $o(h) = o(h').$  
 \end{enumerate} 
Relations of the second type say that $\epsilon(h) h \overline{h}$ 
in the quotient algebra depends only on the outgoing vertex of $h.$ 
We denote $X_a= \epsilon(h) h \overline{h}$ where $a= o(h).$

If the graph has only two vertices, $a$ and $b,$ and one edge connecting 
them, we let $A(\qQ)$ be the quotient of the path algebra by 
relations $ h \overline{h} h=0 =\overline{h} h \overline{h}$ (where 
$h$ is one of the orientations of the edge and $\overline{h}$ is the 
reverse of $h$).  
If the graph has only one vertex $a,$ and no edges, 
define $A(\qQ)$ as the exterior algebra on one 
generator $X_a,$ and place it in degree $2$ to make $A(\qQ)$ graded.
If the graph has more than one vertex, we grade $A(\qQ)$ by lengths 
of paths.
The graded algebra $A(\qQ)$ is (up to isomorphism) independent of the choice
of the orientation $\epsilon$.

For simplicity, we will write $A$ instead of $A(\qQ)$.
The algebra $A$ is finite-dimensional, $\mbox{dim}(A)=2(|I|+|E|).$ Any path 
of length at least 3 equals $0$ in $A.$ 

Notice that $X_a$ (see above) is central, and $X_a,$ over all vertices $a,$ 
form a basis for the degree $2$ subspace of $A.$ 

A length $0$ path $(a),$ for a vertex $a\in I,$ is a minimal idempotent 
in $A,$ and $1=\sum_a (a).$ 

\vspace{.1in}

\noindent
\emph{Example:} If $\qQ$ has one vertex and one loop, $A$ is isomorphic 
to the exterior algebra on two generators $h, \overline{h}$: 
\begin{equation*} 
  A\cong \C \langle h, \overline{h}\rangle / h^2=\overline{h}^2=h\overline{h}+ 
  \overline{h} h = 0. 
\end{equation*} 

The trace $tr: A \to \C$ defined by
$$tr(X_a)=1, \qquad tr(h )=tr((a))=0,$$
makes $A$ into a graded Frobenius algebra. Note that $A$ is symmetric (but with 
respect to a different trace) if and 
only if $\qQ$ is bipartite (i.e if it is possible to
partition the set of vertices of $\qQ$ into two disjoint sets
in such a way that all edges go from one set to the other).
In the latter case, $A$ is isomorphic to the 
zigzag algebra of $\qQ$, see [HK]. For any $\qQ$, the algebra $A$ is 
a skew-zigzag algebra, in the terminology of [HK, Section 4.6].

\paragraph{}Denote by $P_a$ the indecomposable projective left $A$-module $A(a).$ 
An indecomposable projective left $A$-module is isomorphic to $P_a,$ for 
some $a.$ Denote by $_aP$ the indecomposable projective right 
$A$-module $(a)A.$ Since $A$ is Frobenius, $P_a$ and $\hsm _aP$ 
are, in addition, injective. 

\vspace{.1in}

Let $S_a$ be the simple quotient of $P_a$ (equivalently, the quotient 
of $P_a$ by all paths of length greater than $0$).  Denote by 
 $\widehat{a}$ 
the image of $(a) \in P_a$ under the quotient map. $S_a$ is a one-dimensional 
complex vector space and is spanned by $\widehat{a}.$ 
A simple left $A$-module is isomorphic to $S_a,$ for some $a.$ 
The modules $P_a, S_a, $ and $_aP$ inherit $\Z$-grading from $A.$ 

\vspace{.1in}

Denote by $[m]$ the grading shift down by $m$., i.e $(M[m])^l=M^{m+l}$.
Let $\mathbf{Mod}(A)$ (resp. $\mathbf{HMod}(A)$) be the category
of graded $A$-modules (resp. the category of graded $A$-modules
equipped with a Hermitian structure $x \mapsto x^*$ such that
$h({x}^*)=(hx)^*$ for all edges $h$). For any two graded $A$-modules
$M,N$ we denote by $\mathrm{Hom}_A(M,N)$ the set of grading-preserving
$A$-homomorphisms. 

\paragraph{}The modules $P_a, \hspace{0.03in} _aP, S_a$ have unique Hermitian 
structures $x \mapsto x^*$ such that $(a)^*=(a)$ and $h(x^*)=(hx)^*$ for all edges $h$.

\vspace{.1in}

\paragraph{5.2.} Let $V_a, W_a$, $a \in I$ be finite-dimensional $\C$-vector
spaces. Consider the graded $A$-module 
$$M= \bigoplus_a V_a\otimes P_a \oplus W_a\otimes S_a[-1],$$ 
a direct sum of a projective and a semisimple $A$-module. We raise 
the grading of simple modules $S_a$ by $1$ to "balance" them in the 
middle of projective modules $P_a,$ the latter non-zero in degrees $0,1,2.$ 

Let's equip $M$ with a degree $1$ generalized differential 
$d: M\to M$ with the property $d^2=c$ for a fixed degree $2$ central 
element $c$ of $A,$ 
 $$  c = \sum_a c_a X_a, \hspace{0.3in} c_a\in \C. $$ 
 $d$ should super-commute with $A,$ 
  $$ dx = (-1)^{|x|} xd, \hspace{0.1in} x\in A.$$ 
The graded components of $M$ are 
 \begin{eqnarray*} 
   M^0 & = & \oplusop{a} (V_a \otimes (a)), \\
   M^1 & = & \bigoplusop{a} ( (W_a \otimes \widehat{a}) 
   \oplusop{o(h)=a} (V_{i(h)} \otimes h ) ),  \\
   M^2 & = & \oplusop{a} (V_a \otimes X_a). 
 \end{eqnarray*} 
and the differential must have the form
$$ 0 \lra M^0 \stackrel{d^0}{\lra} M^1 \stackrel{d^1}{\lra} M^2 \lra 0. $$ 
Since $d(a)= (a) d,$ for minimal 
idempotents $(a)\in A,$ the generalized complex decomposes into the sum of 
$$ 0 \lra (a)M^0 \stackrel{d^0}{\lra} (a)M^1 \stackrel{d^1}{\lra} (a)M^2 
  \lra 0, $$
over all $a.$ We can write the latter as 
 $$ 0 \lra V_a \otimes (a) \stackrel{d^0}{\lra} 
  \oplusop{o(h)=a} (V_{i(h)} \otimes h )
  \oplus (W_a \otimes \widehat{a}) 
  \stackrel{d^1}{\lra}  V_a \otimes X_a \lra 0. $$ 
 The components of $d^0$ can be described as maps   
 $B_{h}\in \mbox{Hom}(V_a, V_{i(h)}), j_a \in \mbox{Hom}(V_a, W_a)$:
 $$d^0= (\oplusop{o(h)=a} B_{h}, j_a)^t.$$  
 The upper $t$ in the formula stands for  transposing a row vector 
 into a column vector. From $d h = - h d,$ for all edges $h,$ we 
 derive that 
 $$d^1= (\oplusop{o(h)=a} \epsilon(\overline{h}) 
 B_{\overline{h}}, i_a),$$ 
 where $i_a\in \mbox{Hom}(W_a,V_a).$ 
 We should have $c=d^1d^0,$ or, specializing to a vertex $a,$ 
 \begin{equation}\label{centerm}
 c_a\mathrm{Id} = \sum_{o(h)=a} 
  \epsilon(\overline{h}) B_{\overline{h}}B_{h} + i_aj_a.
 \end{equation}  
 The right hand side is the $a$-component of the complex moment map 
 for the Nakajima quiver varieties. 
 
\vspace{0.1in}

If $d$ is given as above by the data $d=(B_h, i_a, j_a)$ we may define its
Hermitian adjoint $d^*=(\epsilon(\overline{h})B^*_{\overline{h}}, -j_a^*, i_a^*)$
to be of the same form. The real component of the moment map equation
$\mu_{\mathbb{R}}(B,i,j)=\boldsymbol{\zeta}_{\R}$ is equivalent to the relation
\begin{equation}\label{E:centerm2}
\frac{\sqrt{-1}}{2}(dd^*+d^*d)=\sum_a \zeta_{\R,a} X_a.
\end{equation}
We can think of $\boldsymbol{\zeta}_{\mathbb{C}}$ and  
$\boldsymbol{\zeta}_{\mathbb{R}}$ as degree 2 central elements of 
 $A$, by taking the standard bases of $\mathbb{C}^I$ and 
$\mathbb{R}^I$ to $\{ X_a\}_{a\in I}$. Collapse the grading
from $\Z$ to $\Z/2\Z$, and write $\mathbf{Mod}_2(A)$ and $\mathbf{HMod}_2(A)$ for the
corresponding categories of $\Z/2\Z$-graded modules. Equations~\ref{centerm} and \ref{E:centerm2}
together with the definitions of quiver varieties imply the following result.

\begin{prop} \label{bijall}
For any $\zeta_{\C}, \zeta_{\RR}$ there is 
a bijection between points on the Nakajima variety 
$\mathcal{M}_{\zeta}(\mathbf{v},\mathbf{w})$ and isomorphism classes 
of the following data $(M,d, \psi)$: 

A graded $A$-module $M \in \mathbf{HMod}_2(A)$ which is a direct sum of a projective and a 
semisimple module, with $v_a$ the multiplicity of $P_a$ and $w_a$ the 
multiplicity of $S_a[-1],$ with a generalized differential $d$ such 
that $d^2=\boldsymbol{\zeta}_{\C},$ 
and $\frac{\sqrt{-1}}{2}(dd^*+d^*d)=\boldsymbol{\zeta}_{\mathbb{R}},$
and isomorphisms $\psi_a:W_a\cong 
\Hom_A(S_a[-1], M)$.
\end{prop} 

Now suppose that $\zeta_{\C}$ is generic. The complex
description of quiver varieties yield the following result.

\begin{prop} \label{bij1} There is 
a bijection between points on the Nakajima variety 
$\mathcal{M}_{\zeta}(\mathbf{v},\mathbf{w})$ and isomorphism classes 
of the following data $(M,d, \psi)$: 

A graded $A$-module $M \in \mathbf{Mod}_2(A)$ which is a direct sum of a projective and a 
semisimple module, with $v_a$ the multiplicity of $P_a$ and $w_a$ the 
multiplicity of $S_a[-1],$ with a generalized differential $d$ such 
that $d^2=\boldsymbol{\zeta}_{\C},$ 
and isomorphisms $\psi_a:W_a\cong 
\Hom_A(S_a[-1], M)$.
\end{prop} 

\paragraph{}It is easy to check that two nonisomorphic data $(M,d,\psi)$ as above
remain nonisomorphic after applying the functor $\Phi_1$ (see (\ref{E:fonct})),
and that 
$$\mathrm{Hom}_A(S_a[-1],M)\cong \mathrm{Hom}_{\UMod_2(A)}
(S_a[-1],R\Phi_1(M))$$
for any $M$ as above. This gives the following variant of Proposition~\ref{bij1}.

\begin{prop}\label{bij22} There is 
a bijection between points on the Nakajima variety 
$\mathcal{M}_{\zeta}(\mathbf{v},\mathbf{w})$ and isomorphism classes 
of the following data $(\underline{M},\psi)$: 

An object $\underline{M} \in \UMod_2(\wA_{\boldsymbol{\zeta}_{\C}})$ such that
$\underline{M}\simeq\Phi_1(M')$ for some $M' \in \mathbf{Com}_2(A,\boldsymbol{\zeta}_{\C})$
with $M'\simeq \bigoplus_s S_a[-1] \otimes \C^{w_a} \oplus \bigoplus_a P_a \otimes \C^{v_a}$
as an $A$-module;
and isomorphisms $\psi_a:W_a\cong 
\Hom_{\UMod_2(A)}(S_a[-1], R(\underline{M}))$.
\end{prop} 

\vspace{.1in}

Consider now the case when $\zeta_\C=0$ and $\zeta_{\R} \in (\N^+)^I.$ The description
of $\mathbf{M}^{ss}_{\zeta}(\mathbf{v},\mathbf{w})$ given in Section~1.3 i)
yields the following. 

\begin{prop} There is a bijection between points on 
$\mathcal{M}_{\zeta}(\mathbf{v},\mathbf{w})$ and isomorphism classes of 
$(M,d,\psi)$ where 
$M \in \mathbf{Mod}_{2}(A)$ and $\psi$ are as in Proposition~\ref{bij1},
$d^2=0$ and no projective submodule of $M$ 
is $d$-stable. 
\end{prop} 

\vspace{0.1in}

\noindent
\textbf{Remark.} All the above results also hold in the $\Z$-graded case. 

\section{Weyl group action in categories of duplexes}

\paragraph{6.1.} We use notations from Section 5.1. Choose a vertex $a$
without a loop. Since $P_a$ is a left $A$-module and $_aP$ a right 
$A$-module, $P_a\otimes \hsm _a P$ is an $A$-bimodule (the tensor 
product is over $\C$).  Denote by $m_a: P_a \otimes \hsm _aP \lra A$ 
the bimodule map which is simply the restriction of the multiplication 
map $A\otimes A \lra A,$ so that $m_a ((a)\otimes (a)) = (a).$
Denote by $\Delta_a : A\lra P_a \otimes \hsm _aP$ the bimodule map
determined by
  \begin{equation*}
    \Delta_a(1) = X_a\otimes (a) + (a) \otimes X_a + \sum_{o(h)=a} 
 \epsilon(h) \overline{h}\otimes h,  
    \end{equation*}
the sum over all edges $h$ that start at $a.$ 

\paragraph{6.2.}
Let $C_{a,x},$ for $x\in \C,$ be the following duplex of bimodules:
\begin{equation}\label{main-du}
  \stackrel{\Delta_a}{\lra} P_a \otimes \hsm _aP
  \stackrel{ x m_a}{\lra} A \stackrel{\Delta_a}{\lra}
\end{equation}

If $x\not= 0,$ this duplex is isomorphic to
\begin{equation}\label{main2-du}
  \stackrel{x \Delta_a}{\lra} P_a \otimes \hsm _aP
  \stackrel{ m_a}{\lra} A \stackrel{x \Delta_a}{\lra}
\end{equation}

Denote by $d$ the differential in (\ref{main-du}) and (\ref{main2-du}). We have
  \begin{eqnarray*}
  d^2\mid_{P_a\otimes\hsm _aP} & = & x (X_a \otimes 1 + 1\otimes X_a), 
     \\
  d^2\mid_A                    & = & x (2 X_a - \sum_{o(h)=a} X_{i(h)}).
  \end{eqnarray*}
Here $X_a\otimes 1 $ is the left multiplication by $X_a,$ etc...
 Hence, as operators on $C_{a,x}$ we have
  \begin{equation} \label{d1-square}
   d^2 = x ( X_a\otimes 1 + 1 \otimes X_a - \sum_{o(h)=a} X_{i(h)}
 \otimes 1  ) 
  \end{equation}
Note that $X_b \otimes 1 - 1 \otimes X_b$ acts trivially on $P_a \otimes \hsm _aP$
and $A$ if $a \neq b$. Thus we also have
\begin{equation} \label{d-square}
   d^2 = x ( X_a\otimes 1 + 1 \otimes X_a - \sum_{o(h)=a} X_{i(h)}
 \otimes 1  ) 
 + \sum_{b\not= a} x_b (X_b\otimes 1 - 1 \otimes X_b)
\end{equation}
for any $x_b\in \C,$ as $b$ ranges over all vertices other than $a.$ 

\vspace{0.1in}

The Weyl group $\mathbf{W}$ of $(I,E)$ has generators $s_a,$ over 
all loopless vertices $a\in I,$ and relations  

 (i) $s_a^2=1$ for all $a$,  

 (ii) $s_a s_b = s_a s_b$ if $a$ and $b$ do not have a common edge, 

 (iii) $s_a s_b s_a = s_b s_a s_b$ if $a$ and $b$ are joined by exactly 
 one edge. 

\vspace{.1in}

The Weyl group acts on $Z_2(A),$ the degree 2 summand of the center of 
$A$ (on the vector space with the basis $\{ X_a\}_{a\in I}$) by 
 \begin{equation*} 
  s_a(c) = c + x_a (\sum_{o(h)=a}X_{i(h)} - 2 X_a)
 \end{equation*} 
for $c=\sum_{b\in I} x_b X_b.$ 
This action is compatible with the one defined on
$\zeta_{\C}$ in Section~1.6. via the natural identification
$\C^I \to Z_2(A)$.
It follows from (\ref{d-square}) 
that, for any $c\in Z_2(A)$, the tensor product with $C_{a,-x_a}$ is a 
functor from $\mathbf{Com}_2(A,c)$ to $\mathbf{Com}_2(A,s_a(c))$ and 
from $\mc{K}_2(A,c)$ to $\mc{K}_2(A,s_a(c))$. Denote this functor 
by $\mc{R}_a:\mc{K}_2(A,c) \to \mc{K}_2(A,s_a(c)).$

\begin{lem} The functor $\mathcal{R}_a$ lifts to a functor
from $\UMod_2(\wA_c)$ to $\UMod_2(\wA_{s_a(c)})$.
\end{lem}
\noindent
\textit{Proof.} Let $\tilde{P}_b=\wA_c \otimes_A P_b$ be an indecomposable projective
$\wA_c$-module. We have to show that $C_{a,-x_a} \otimes \tilde{P}_b$ is a projective
$\wA_{s_a(c)}$-module. By definition, we have
$$C_{a,x} \otimes \tilde{P}_b= \begin{pmatrix} A \otimes \tilde{P}_b \\ \oplus \\ 
(P_a \otimes \hsm _a P)[1] \otimes_A
\tilde{P}_b \end{pmatrix}=\begin{pmatrix} P_b \oplus P_b[1]\\
\oplus\\\bigoplus_{\underset{o(h)=a}{i(h)=b}}(P_a[1] \oplus P_a) \end{pmatrix}.
$$
and the action of the element $d \in A_{s_a(c)}$ is given by
$$d=\begin{pmatrix} 1 \otimes d & m_a \otimes 1 \\
\Delta_a \otimes 1 & 1 \otimes d \end{pmatrix}.$$
Since $A$ is Frobenius, $P_a$ and $P_b$ are both projective and injective.
Furthermore, it is easy to check that $d$ is injective on $P_b$ and that
$P_b \cap d(P_b)=0$. Similarly, we have $d$ is injective on $P_a[1]$ and
we have $P_a[1] \cap d(P_a[1])=0$. By applying Proposition~\ref{quot-by-projinj}
twice we see that $C_{a,x} \otimes \tilde{P}_b=0$ in $\UMod_2(\wA_{s_a(c)})$ as
desired.\qed

\vspace{.1in}

We will denote this functor by the same symbol 
$\mathcal{R}_a:\UMod_2(\wA_c)\to\UMod_2(\wA_{s_a(c)})$.

\vspace{.1in}

\paragraph{6.3.} We now deal with the the categorical analogue of the
braid relation.
\begin{prop} \label{inv-pr} 
 If $x \not= 0,$ there is an isomorphism in the stable
category of bimodule duplexes
\begin{equation} \label{invertible}
C_{a,-x}\otimes_A C_{a,x} \cong A.
\end{equation}
\end{prop}
\noindent
\emph{Proof:} Let $N= C_{a,-x}\otimes_A C_{a,x}$ and $\partial$
be the differential in $N.$ Notice that $\partial^2=0,$ and $N$ is a
duplex of $A$-bimodules.

Since $x\not= 0,$ we can write $C_{a,-x}$ as
\begin{equation}
   \stackrel{\Delta_a}{\lra} P_a \otimes \hsm _aP
   \stackrel{- x  m_a}{\lra} A \stackrel{\Delta_a}{\lra}
\end{equation}

$N$
is the total duplex of the following diagram (which is 2-cyclic in
horizontal and vertical directions, and each of the four squares anticommutes)

 \begin{figure}[htb] \label{double}
 \begin{pspicture}(-1,0)(6,6)
 \begin{psmatrix}
 \hsm  &  \hsm     &  \hsm     &     \hsm    \\
 \hsm  & $A\otimes_A P_a\otimes \hsm _aP$   &  $ P_a\otimes \hsm _aP\otimes_A
     P_a\otimes \hsm _aP$     &     \hsm    \\
 \hsm  &  $A\otimes_A  A$     &  $P_a\otimes \hsm _aP\otimes_A  A$ 
&  \hsm  \\
 \hsm  &   \hsm    &   \hsm    & \hsm
 \psset{arrows=->}
   \ncline{2,1}{2,2}^{\scriptsize{$-x m_a \otimes \id\hspace{0.5in}$}}
   \ncline{2,2}{2,3}^{\scriptsize{$\id\otimes \Delta_a\hspace{0.2in}$}}
   \ncline{2,3}{2,4}^{\scriptsize{$\hspace{0.8in}-x m_a \otimes \id$}}
   \ncline{3,1}{3,2}^{\scriptsize{$-x m_a \otimes \id\hspace{0.3in}$}}
   \ncline{3,2}{3,3}^{\scriptsize{$\id\otimes \Delta_a$}}
   \ncline{3,3}{3,4}^{\scriptsize{$\hspace{0.4in}-x m_a \otimes \id$}}
   \ncline{1,2}{2,2}>{\scriptsize{$-\id \otimes x\Delta_a$}}
   \ncline{2,2}{3,2}>{\scriptsize{$-\id\otimes m_a$}}
   \ncline{3,2}{4,2}>{\scriptsize{$-\id \otimes x \Delta_a$}}
   \ncline{1,3}{2,3}>{\scriptsize{$\id \otimes  x \Delta_a$}}
   \ncline{2,3}{3,3}>{\scriptsize{$\id\otimes m_a$}}
   \ncline{3,3}{4,3}>{\scriptsize{$\id \otimes  x \Delta_a$}}
 \end{psmatrix}
  \end{pspicture}
  \end{figure}

\begin{figure} [htb] \drawing{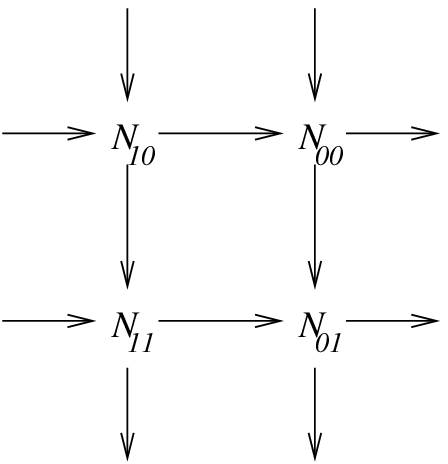}\caption{} \label{dupn} 
 \end{figure}

Denote by $N_{ij}$ the four bimodules in the above diagram, see 
 figure~\ref{dupn}. 

Simplifying our notation as at the end of Section 3.1, we 
write $N$ as
  \begin{equation*}
   \stackrel{\partial}{\lra} N_{00}\oplus N_{11} 
   \stackrel{\partial}{\lra}
   N_{01}\oplus N_{10} \stackrel{\partial}{\lra}. 
  \end{equation*}

Notice that $_aP\otimes_A P_a\cong \C(a) \oplus \C  X_a$ is a 
two-dimensional vector space. Let
$\zeta: N_{11}\to N_{00}$ be the map
\begin{equation*}
N_{11}\cong A \stackrel{-\Delta_a}{\lra} P_a\otimes \hsm _aP \lra
P_a\otimes \hspace{0.03in} _aP\otimes_A P_a\otimes \hspace{0.03in} _aP,
\end{equation*}
where the last map takes $u_1\otimes u_2$ to $u_1\otimes (a) 
 \otimes u_2.$ Let $N'= \{ u + \zeta(u)| u\in N_{11} \}.$ It is a 
 subbimodule of $N_{11}$ isomorphic to $A,$ and $\partial N'=0.$

Let $N''=P_a \otimes (a) \otimes \hsm _aP.$ It is a subbimodule of 
 $N_{00}.$

\begin{lem} $N$ is a direct sum of its 3 subduplexes
  \begin{eqnarray*}
    T_{-1} & = & \{  \lra N_{10} \lra \partial N_{10} \lra  \}, \\
    T_0    & = & \{  \lra N' \lra 0 \lra  \}, \\
    T_1    & = & \{  \lra N''  \lra \partial N'' \lra \}.
  \end{eqnarray*}
\end{lem}

Since $\partial$ is injective on $N_{10}$ and on $N''$, and both
$N_{10}$ and $N''$ are projective bimodules, the duplexes
$T_{-1}$ and $T_1$ are stably equivalent to the zero duplex.
Therefore, $N$ is equivalent in the stable category to the bimodule 
 duplex $\{ \lra A\lra 0 \lra \}.$ \qed

\vspace{0.1in}

\noindent
\textbf{Remarks} i) Proposition~\ref{inv-pr} says that the functor $\mc{R}^2_a:$ 
 \begin{equation*} 
  \mc{K}_2(A,c) \stackrel{\mc{R}_a}{\lra} \mc{K}_2(A,s_a(c)) 
  \stackrel{\mc{R}_a}{\lra} \mc{K}_2 (A,c).  
 \end{equation*} 
is isomorphic to the identity functor, as long as 
 $s_a(c)\not= c$ (equivalently, if $c_a\not= 0).$ \\
ii) The isomorphism (\ref{invertible}) holds for $x=0$ as well, if 
we use (\ref{main-du}), with $x=0,$ to define one of the duplexes on 
the left hand side of (\ref{invertible}) and (\ref{main2-du}) to 
define the other.   

\begin{prop} If $a$ and $b$ are not connected by an edge, for any 
$x,y\in \C$ there is an isomorphism of bimodule duplexes
\begin{equation}\label{easy-iso}
    C_{a,x } \otimes_A C_{b, y } \cong C_{b, y }\otimes_A C_{a, x}.
\end{equation}
\end{prop}
\noindent
\emph{Proof:} Since $_aP \otimes_A P_b\cong 0 \cong
  \hsm  _b P \otimes_A P_a,$ left and right hand sides of 
 (\ref{easy-iso}) are isomorphic to
  \begin{equation*}
  \begin{CD}
  @>{m_a + m_b}>>  A  @>{x \Delta_a + y \Delta_b}>>
    (P_a\otimes \hsm _a P) \oplus (P_b\otimes \hsm _b P)
    @>{m_a + m_b}>>
  \end{CD}
   \end{equation*}
\qed

\begin{prop} If $a$ and $b$ are connected by one edge, for any 
 $x,y\in \C$ 
there is an isomorphism in the stable category of bimodule duplexes
\begin{equation}\label{hard-eq}
  C_{a,y}\otimes_A C_{b,x+y} \otimes_A C_{a,x} \cong
  C_{b,x}\otimes_A C_{a,x+y} \otimes_A C_{b,y}.
\end{equation}
\end{prop}
\noindent
\emph{Proof:} Denote by $N$ the duplex on the left hand side of
(\ref{hard-eq}). It is built out of 8 bimodules
\begin{equation*}
   N_{ijk} = C_{a,y}^i\otimes_A C_{b,x+y}^j \otimes_A
   C_{a,x}^k, \hspace{0.2in} i,j,k\in \{ 0,1\}.
\end{equation*}

The differential $\partial$ of $N$ is injective on $N_{000}$ (already 
the component $N_{000}\to N_{010}$ of $\partial$ is injective). Let
  \begin{equation*}
    T_{-1} = \{ \stackrel{\partial}{\lra} N_{000}
    \stackrel{\partial}{\lra} \partial N_{000} 
    \stackrel{\partial}{\lra} \}
  \end{equation*}
be the subduplex of $N$ generated by $N_{000}.$ 

Let $N'= P_a\otimes (a) \otimes \hsm _a P\subset P_a \otimes \hsm _a P
   \otimes_A P_a \otimes \hsm _a P \cong N_{010}.$ The differential 
is injective on $N_{010}$ (since the component $N_{010}\to N_{110}$ of
  $\partial$ is injective). Let
\begin{equation*}
   T_1 = \{ \stackrel{\partial}{\lra} \partial N'
    \stackrel{\partial}{\lra} N' \stackrel{\partial}{\lra} \}
  \end{equation*}
be the subduplex of $N$ generated by $N_{010}.$ 

The algebra $A$ is Frobenius, and each projective $A$-module is 
 injective. In particular, $P_a,P_b$ are injective $A$-modules, and
 $N', N_{000}$ (both isomorphic to $P_a\otimes\hsm _a P$) are injective
 $A\otimes A^o$-modules (that is, injective $A$-bimodules). Moreover,
 $T_1 \cap T_{-1}= 0.$ Applying proposition~\ref{quot-by-projinj}
 twice, we see that duplexes $N$ and $\widetilde{N}= N/(T_{-1}\oplus 
  T_1)$
  are isomorphic in the stable category of duplexes of bimodules.

Let $h$ be the edge with $o(h)=a$ and $i(h)=b.$ 
    The duplex $\widetilde{N}$ is isomorphic to

  \begin{equation} \label{duplex-N}
  \stackrel{\widetilde{\partial}^1}{\lra}
   \left( \begin{array}{c}
     P_a \otimes \hsm _a P \\
     \oplus                \\
     P_b \otimes  \hsm _b P
   \end{array} \right)
   \stackrel{\widetilde{\partial}^0}{\lra}
   \left( \begin{array}{c}
      P_a \otimes \hsm _b P \\
      \oplus                \\
      A                     \\
      \oplus                \\
       P_b \otimes  \hsm _a P
   \end{array} \right)
   \stackrel{\widetilde{\partial}^1}{\lra}
   \end{equation}
with the differential $\widetilde{\partial}$ given by matrices
  of bimodule maps

  \begin{equation} \label{first-matrix}
    \widetilde{\partial}^0= \left( \begin{array}{cc}
   \epsilon(\overline{h})y \hspace{0.05in}\id\otimes \overline{h} &  
   \epsilon(h) y \overline{h}\otimes \id    \\
                        m_a      &    -m_b      \\
   \epsilon(h) x h \otimes\id  &   \epsilon(\overline{h})x 
    \hspace{0.05in}\id\otimes h 
                           \end{array}\right)
  \end{equation}

   \begin{equation} \label{last-matrix}
     \widetilde{\partial}^1= \left( \begin{array}{ccc}
     \id \otimes h    &  (x+y)\Delta_a   &  -\overline{h}\otimes \id  \\
     h \otimes \id     &  -(x + y)\Delta_b  &  -\id \otimes \overline{h}
                                     \end{array}\right)
   \end{equation}

The following example explains our notations: the entry
 $\epsilon(\overline{h})y \hspace{0.05in}\id\otimes \overline{h}$ in 
 the top left corner of (\ref{first-matrix}) is a bimodule map 
 $P_a\otimes \hsm _a P \lra P_a\otimes \hsm _b P$ which takes 
 $u\otimes v\in P_a\otimes \hsm _a P$ to $\epsilon(\overline{h})y 
 u\otimes \overline{h}v\in P_a\otimes \hsm _b P.$

Denote by $M$ be the duplex on the right hand side of (\ref{hard-eq}).
Since the right hand side is obtained from the left hand side by
switching $a$ with $b,$ $x$ with $y,$ and $h$ with $\overline{h}$ we 
see that $M$ is isomorphic in the stable category to the duplex 
$\widetilde{M}$ defined by making these switchings  
in formulas (\ref{duplex-N}), (\ref{first-matrix}), (\ref{last-matrix}).
It is easy to check that duplexes $\widetilde{N}$ and $\widetilde{M}$
are isomorphic. Therefore, duplexes $N$ and $M$ are isomorphic in
the stable category of duplexes. \qed

\vspace{.1in}

 We may restate the above results in the following form.

\begin{theorem}\label{T:2} The functors $\mc{R}_a: \UMod_2(\wA_c) \to \UMod_2(\wA_{s_a(c)})$ define a braid group action 
on the family of categories $\UMod_2(\wA_{w(c)})$ for $w \in \mathbf{W}$; in other words, we have, for any
$c \in Z_2(A)$ isomorphisms of functors
$$\mathcal{R}_a \circ \mathcal{R}_b \circ \mathcal{R}_a \simeq \mathcal{R}_b \circ
\mathcal{R}_a \circ \mathcal{R}_b: \UMod_2(\wA_c) \to \UMod_2(\wA_{s_as_bs_a(c)})$$
if $a$ and $b$ are connected by one edge, and
$$\mathcal{R}_a \circ \mathcal{R}_b \simeq \mathcal{R}_b \circ
\mathcal{R}_a :  \UMod_2(\wA_c) \to \UMod_2(\wA_{s_as_b(c)})$$
if $a$ and $b$ are not connected.
Moreover, if $c$ lies in a generic orbit then this action factors through a Weyl group
action, i.e we have
$$\mathcal{R}_a \circ \mathcal{R}_a \simeq Id: \UMod_2(\wA_c) \to \UMod_2(\wA_c)$$
for any $a \in I$.
\end{theorem}

Passing to the homotopy categories $\mathcal{K}_2(A,c)$ we obtain
 
\begin{theorem} The functors $\mc{R}_a: \mathcal{K}_2(A,c) \to \mathcal{K}_2(A,s_a(c))$ define a braid group action 
on the family of categories $\mathcal{K}_2(A,w(c))$ for $w \in \mathbf{W}$. This action factors through to a Weyl group
action if $c$ lies in a generic orbit.\end{theorem}

\vspace{.1in}

\noindent
\textbf{Remark.} In order to ensure that
the above braid group action on the set of categories $\UMod_2(\wA_{w(c)})$
or $\mathcal{K}_2(A,w(c))$ factors through the Weyl group it is enough to
impose the following weaker condition : no point of the orbit of $c$ under $W$
is fixed by any of the reflexions $s_a$ for loopless vertices $a$.

\section{Weyl group actions on Nakajima varieties}

\paragraph{7.1.} Let us say that 
$\boldsymbol{\zeta}_\C=-\sum_a\zeta_{\C,a} X_a \in A$ is \textit{generic} if 
$\zeta_\C$  is generic in the
sense of Section~1.3.  If $N$ is any $(A,\boldsymbol{\zeta}_{\C})$-duplex with 
$\boldsymbol{\zeta}_\C$ generic
and $a \in I$ we set $\mc{R}_a(N)= C_{a,\zeta_{\C,_a}}\otimes N$.
%where $\zeta_a=-\langle \alpha_a,\lambda\rangle$.

It follows from the results in Section~6 that this defines
an action of the Weyl group $\mathbf{W}$ on the set of objects of
$\UMod_2(\wA_c)$ \textit{for all generic} $c$.

\paragraph{}So let us assume that $\zeta_{\C}$ is generic,
and let us identify the points of
$\mathcal{M}_{\zeta}
(\mathbf{v},\mathbf{w})$ with isomorphism classes of data $(\underline{M},\psi)$
as in Proposition~\ref{bij22}.

\begin{theorem}\label{T:7.1} The functor $\mc{R}_a$ induces a bijection 
of sets from $\mathcal{M}_{\zeta}
(\mathbf{v},\mathbf{w})$ to $\mathcal{M}_{s_a(\zeta)} (s_a(\mathbf{v},\mathbf{w}))$, 
which coincides with the isomorphism $\kappa_{s_a}$.\end{theorem}
\noindent
\textit{Proof.} Let us describe the action of $\mc{R}_a$ in more details. 
 For notational convenience we will write $N^{\oplus V}$ for the 
tensor product $N \otimes_\C V$ when $N$ is an $\wA_c$-module and 
$V$ a $\C$-vector space. We will also denote by $(M,d)$ the objects of
$\UMod_2(\wA_c)$. If $(M,d) \in \UMod_{2}(\wA_{\underline{\zeta_{\C}}})$ 
then by definition
$\mc{R}_a(M,d)=(M',d') \in \UMod_{2}(\wA_{s_a\underline{\zeta_{\C}}})$
where
$$M'=\begin{pmatrix} A \otimes M \\ \oplus \\ (P_a \otimes \hsm _a 
P)[1] \otimes_A
M \end{pmatrix}$$
and
$$d'=\begin{pmatrix} 1 \otimes d & m_a \otimes 1 \\
\Delta_a \otimes 1 & 1 \otimes d \end{pmatrix}.$$
Let us write
$$M=\bigoplus_{k \in I} (P_k^{\oplus V_k} \oplus S_k^{\oplus W_k}[1]).$$
Observe that
$$A \otimes P_a^{\oplus V_a} \stackrel{\Delta_a \otimes 1}{\lra} 
(P_a \otimes \hsm_a P)[1] \otimes_A P_a^{\oplus V_a}$$
is injective. This implies that $A \otimes P_a^{\oplus V_a} \oplus
d'(A \otimes P_a^{\oplus V_a})$ is stably trivial. Hence, by Proposition~\ref{quot-by-projinj}, 
$(M',d') \simeq (M'',d'')$ where $M''=M'/(A \otimes P_a^{\oplus V_a} 
 \oplus  d'(A \otimes P_a^{\oplus V_a}))$.
A direct computation shows that
\begin{equation}\label{E:7.1}
(P_a \otimes \hsm _a P)[1] \otimes_A M=
\begin{pmatrix}
P_a^{\oplus V_a}[1]\\
\oplus\\
P_a^{\oplus V_a}[1]\\
\oplus\\
P_a^{\oplus W_a}\\
\oplus\\
\underset{b-a}{\bigoplus} P_a^{\oplus V_b}
\end{pmatrix}
\end{equation} and
$$M''\simeq \begin{pmatrix}
\underset{k\neq a}{\bigoplus} P_k^{\oplus V_k} \\
\oplus\\
\underset{k\neq a}{\bigoplus} S_k^{\oplus W_k}\\
\oplus\\
S_a^{\oplus W_a}\\
\oplus\\
P_a[1]^{\oplus V_a}\\
\oplus\\
  P_a^{\oplus W_a}\\
\oplus\\
\underset{b-a}{\bigoplus} P_a^{\oplus V_b}
\end{pmatrix}$$
We have
\begin{equation} \label{big-matrix}     d'' =
\begin{pmatrix}
B & i & 0 & 0 & -Bi & \lambda_a Id -\epsilon BB\\
j & 0 & 0 & 0 & 0 & 0\\
0 & 0 & 0 & 0 & \lambda_a Id- ji & -jB\\
0 & 0 & 0 & 0 & - i & 0\\
0 & 0 & Id & j & 0 & 0\\
Id & 0 & 0 & B & 0 & 0
\end{pmatrix}.
\end{equation}
 From the relation $\mu(B,i,j)=\lambda$ and from the fact that 
 $\lambda$ is generic it follows that the fourth column of 
 (\ref{big-matrix}) is nonsingular.  In particular, 
 $d''_{|P_a[1]^{\oplus V_a}}$ is injective and $P_a[1]^{\oplus V_a}
\oplus d''(P_a[1]^{\oplus V_a})$ is stably trivial.

\vspace{0.2in} 

Thus $(M'',d'')$ is isomorphic, in $\UMod_{2}(\wA_{s_a(\underline{\zeta_{\C}})})$, to
\begin{equation}\label{E:7.9}
(M''',d''')=
(M''/\left( P_a[1]^{\oplus V_a } \oplus d''(P_a[1]^{\oplus V_a}),d'')
\right).
\end{equation}
Moreover, there is a canonical isomorphism $u: R((M,d)) \simeq R((M',d')) \simeq R((M''',d'''))$ and
we may set $\psi'=u \circ \psi$.
Let
$$(B',i',j') \in \mathcal{M}_{s_a(\zeta)} (s_a(\mathbf{v},\mathbf{w}))$$ 
be the element corresponding to $(M''', d''', \psi')$. Comparing with the
construction of Section~1.6 we see that 
$\kappa_{s_a}(B,i,j)=(B',i',j')$, which proves the Theorem. \qed

\section{Nakajima varieties for affine quivers}

\paragraph{}In this section we restrict ourselves to the case when 
$(I,E)$ is an affine bipartite Dynkin diagram, and reinterpret the above 
construction in terms of the McKay correspondence. This section, and the following
one, can be read independently of the rest of the paper, with the exception of
Section~9.3.
 
\paragraph{8.1.} Let $\{\pm 1\} \subset \G \subset SL(2,\C)$ be a finite group and let 
$\{\rho_a\}_{a \in I}$ be the set of  its irreducible
representations. We also let $\rho_0$ and $\rho$ be the trivial (resp. the
natural $2$-dimensional) representation. Let $\qQ=(I,E)$ be the (unoriented)
affine quiver associated to $\G$ via the McKay correspondence, with $I$ as the set of
vertices and with $T_{ab}$ arrows between $a$ and $b$, where 
$$T_{ab}=\mathrm{dim\;Hom}_\G(\rho_a \otimes \rho, \rho_b).$$

\paragraph{8.2} Let us consider the algebra $A_\Gamma:=\Lambda \C^2 \rtimes
\C[\Gamma]$, and set $\widetilde{A}_\Gamma=
A_\Gamma \sotimes \C[d]/d^2 \simeq \Lambda \C^3 \rtimes
\Gamma$ with relations $dz=-zd$ for any $z \in \C^2$ and
$d\gamma=\gamma d$ for any $\gamma \in \Gamma$.
Both $A_\Gamma$ and
$\widetilde{A}_\Gamma$ are naturally $\Z$-graded. We denote by
$\Mod(A_\Gamma)$ and $\Mod_2(A_\Gamma)$ (resp. $\Mod(\widetilde{A}_\G)$ 
and $\Mod_2(A_\G)$) the categories of 
$\Z$-graded and $\Z/2\Z$-graded $A_\G$-modules (resp. 
$\widetilde{A}_\G$-modules).
We will reformulate the definition of $\mathcal{M}_{\zeta}
(\mathbf{v},\mathbf{w})$ entirely in terms of
representation theory of $A_\Gamma$ and $\widetilde{A}_\Gamma$.

\paragraph{}The link with the setting of Section~5 is as follows.
Let $I^{\pm}$ be the set
of indices $a$ such that $\rho_a(-1)=\pm 1$. Then $I=I^+ \sqcup I^-$ and
\begin{equation}\label{E:1}
T_{ab} \neq 0 \Rightarrow a \in I^+, b \in I^-\; \mathrm{or}\; a \in I^-,
b\in I^+.
\end{equation}
In particular, the Dynkin diagram $(I,E)$ is bipartite.
Write $A(\qQ)$ for the zigzag algebra corresponding to $(I,E)$.
Recall that $\{\rho_a\}_{a \in I}$ denotes the set of simple
left $\Gamma$-modules. Then $\{\rho_a^*\}$ is the set of \textit{right} 
simple
$\Gamma$-modules. Consider the right projective $A_\Gamma$-modules
$\hsm _a \mathbf{P}=\rho_a^* \otimes \Lambda \C^2$, and put
$\mathbf{P}=\bigoplus_{a} \hsm _a \mathbf{P}$. \\
\hbox to1em{\hfill}It is easy to check that $A(\Gamma) \simeq 
\mathrm{End}_{A_\Gamma}(\mathbf{P})$. Moreover, the functor 
$\mathbf{P} \otimes -$ induces a Morita equivalence
\begin{equation}\label{E:Morita}
\mathbf{Mod}(A_\Gamma) \simeq \mathbf{Mod}(A(\qQ)).
\end{equation}

\paragraph{8.3.}Note 
that $A_\Gamma$ and $\widetilde{A}_\Gamma$ are symmetric
algebras.
In particular, $A_\Gamma$ and $\widetilde{A}_\Gamma$ are self-injective
algebras
(i.e projective and injective objects coincide). If $M$ is a graded 
$A_\Gamma$
(resp. $\widetilde{A}_\Gamma$)-module then the graded dual space $M^*$ is
again an
$A_\Gamma$ (resp. $\widetilde{A}_\Gamma$)-module.

\paragraph{} If $U$ is any $\Gamma$-module, we will
regard $\Lambda \C^2 \otimes U$ and $(\Lambda \C^2 \sotimes \C[d]/d^2)
\otimes U$ as graded
$A_\Gamma$-module, and $\widetilde{A}_\Gamma$-module respectively,
where $\Lambda^0 \C^2 \otimes U$
is placed in degree $0$.
Note that any projective indecomposable $A_\Gamma$-module
(resp. $\widetilde{A}_\Gamma$-module) is isomorphic to $(\Lambda \C^2 \otimes
\rho_a)[n]$  (resp.
$((\Lambda \C^2 \sotimes \C[d]/d^2)\otimes \rho_a)[n]$)
for some $a\in I$ and $n \in \Z$.

\paragraph{} Let us fix a basis $\{x,y\}$ in $\C^2$.
For any $a \in I$, let us fix intertwiners
\begin{equation}\label{E:1.5}
\bigoplus_{(a,b) \in H} \varphi_{ab} :\;\bigoplus_{(a,b) \in H} \rho_b
\stackrel{\sim}{\to} \C^2 \otimes \rho_a.
\end{equation}
Define a function $\epsilon: H \to \C^*$ by the following condition :
$\pi \circ \varphi_{ba}\circ \varphi_{ab} =\epsilon_{(a,b)} x \wedge y$, where
$\pi\;: \C^2 \otimes \C^2 \otimes \rho_b \to \Lambda^2 \C^2 \otimes 
\rho_b$ is the
projection. Note that $\epsilon(h)+\epsilon(\overline{h})=0$ for any $h \in H$.

\paragraph{8.4.} Let $\UMod (A_\Gamma) $ denote the stable category of
$\Z$-graded $A_\Gamma$-modules (see Section~4.).
\begin{lem}\label{L:4.2}
Let $U$ be a $\Gamma$-module and let us consider it as a $A_\Gamma$
module by trivially extending the action to $A_\Gamma$. If $M \simeq U$
in $\UMod (A_\Gamma)$ then $M \simeq U \oplus P$ in
$\mathbf{Mod} (A_\Gamma)$ for some projective module $P$.\end{lem}
\noindent
\textit{Proof.} Let $f:\;M \to U$ and $f':\;U \to M$ such that $ff'
=Id_U$ and $f'f=Id_M$ in $\UMod(A_\Gamma)$. Note that, for any
projective module $P$, any composition of morphisms $U \to P \to U$ is zero.
Hence $\mathrm{Hom}_{\UMod(A_\Gamma)}(U,U)=\mathrm{Hom}_{A_\Gamma}(U,U)$
and $ff'=Id_U$ in $\mathbf{Mod}(A_\Gamma)$. But then $M \simeq U \oplus
\mathrm{Ker}\;f$ in $\mathbf{Mod}(A_\Gamma)$ and the lemma follows.\qed

\paragraph{} Similarly, we let
$\UMod(\widetilde{A}_\Gamma)$
stand for the stable category of $\Z$-graded
$\widetilde{A}_\Gamma$-modules. Replacing $\Z$ by $\Z/2\Z$, we also 
define the categories 
$\UMod_2(A_\Gamma)$ and $\UMod_2(\widetilde{A}_\Gamma)$.
The stable categories
$\UMod(A_\Gamma)$ and $\UMod(\widetilde{A}_\Gamma)$
are endowed with structures of
triangulated categories (see \cite{Happel}).
 
\paragraph{}Note that $\widetilde{A}_\Gamma$ is a free $A_\Gamma$-module. This
gives rise to functors
$$R:\;\UMod(\widetilde{A}_\Gamma) \to \UMod(A_\Gamma),\qquad
R:\;\UMod_2(\widetilde{A}_\Gamma) \to \UMod_2(A_\Gamma).$$

\paragraph{8.5.} In this section we give the realization
of $\mathcal{M}_{\zeta}(\mathbf{v},\mathbf{w})$
in the case 1.3. ii), i.e $\zeta_{\R}$ is arbitrary
 and $\zeta_{\C}$ is generic.\\
\hbox to1em{\hfill}For every $a \in I$ we let $p_a \in \C[\Gamma]$ be the
(central) primitive idempotent corresponding to $\rho_a$. Set
$c_a=x \wedge y \cdot p_a$. Then $\{1\} \cup \{c_a\}_{a \in I}$ forms a basis
of the center of $A_\Gamma$. Consider the following deformation of
$\widetilde{A}_\Gamma$ :
$$\widetilde{A}_{\Gamma,\zeta_{\C}}=(\Lambda \C^2 \rtimes \Gamma) \sotimes
\C[d]/<d^2
-\sum_a \zeta_{\C,a} c_a>.$$
Let
$\UMod_{2}(\widetilde{A}_{\Gamma,\zeta_{\C}})$ be the stable
categories of $\Z/2\Z$-graded
$\widetilde{A}_{\Gamma,\zeta_{\C}}$-modules. As in Section~8.4, the
embedding
$A_\Gamma \subset \widetilde{A}_{\Gamma,\zeta_{\C}}$ gives rise to a restriction
functor $R:\UMod_{2}(\widetilde{A}_{\Gamma,\zeta_{\C}})\to
\UMod_{2}(A_\Gamma)$. As in the undeformed case,
the algebra $\widetilde{A}_{\Gamma,\zeta_{\C}}$ is symmetric and self-injective.

\paragraph{} Let us fix $\mathbf{w} \in \N^I$, $W=\bigoplus_a W_a$ such that
$\mathrm{dim}\;W=\mathbf{w}$ and let
$\mathbb{W}=\bigoplus_a W_a \otimes \rho_a$ be the corresponding
$\Gamma$-module. We will regard $\W$ as a graded $A_\Gamma$-module, where
$\Lambda\C^2$ acts trivially, placed in degree $0$.
Let $\mathcal{N}_\zeta(\mathbf{w})$ be the set of pairs
$(u,M)$ where $M \in \UMod_{2}(\widetilde{A}_{\Gamma,\zeta_{\C}})$
and $u: \W \to R(M)$ is an isomorphism.

\begin{lem}\label{L:4.61} Let $M$ be any
$\widetilde{A}_{\Gamma,\zeta_{\C}}$-module such that
$R(M) \simeq 0$. Then $M$ is projective.\end{lem}
\noindent
\textit{Proof.} By Lemma~\ref{L:4.2} we have $M \simeq (\Lambda \C^2 \otimes \mathbb{V}_1) \oplus
(\Lambda \C^2 \otimes \mathbb{V}_2)[-1]$ for some $\Gamma$-modules $\mathbb{V}_1$ and $\mathbb{V}_2$.
We may assume that 
\begin{equation}\label{E:gargs}
d(\mathbb{V}_1) \subset \C^2 \otimes \mathbb{V}_1.
\end{equation}
Indeed, any $v \in \mathbb{V}_1$ not satisfying (\ref{E:gargs}) generates
a projective submodule $M'$ of $M$.
Consider the linear map $s: \Lambda^2 \C^2 \to \C,\; x\wedge y \mapsto 1$. From
(\ref{E:gargs}) and the relation $dz=-zd$ for $z \in \C^2$ we deduce that
$Tr_{|V}(s \circ d^2)=0$. But by definition,
$$Tr_{|V}(s \circ d^2)=\sum_a \zeta_{\C,a} \mathrm{dim}\;\mathrm{Hom}(\rho_a, \mathbb{V}_1),$$
and the genericity of $\zeta_{\C}$ implies that $\mathbb{V}_1=0$. Similarly,
$\mathbb{V}_2=0$ and we are done.\qed

\begin{theorem}\label{P:4.61}
There is a natural bijection between the set of isomorphism
classes of elements in $\mathcal{N}_\zeta(\mathbf{w})$ and the set of
points of
$\sqcup_\mathbf{v} \mathcal{M}_{\zeta}(\mathbf{v},
\mathbf{w})$.\end{theorem}
\noindent
\textit{Proof.} Let $(u,M) \in \mathcal{N}_\zeta(\mathbf{w})$.
By Lemma~\ref{L:4.2},
we may assume that $M =(\Lambda \C^2 \otimes \V_1[1]) \oplus \W \oplus (\Lambda \C^2
\otimes \V_0)$ as an $A_{\Gamma}$-module. We can also assume that
\begin{equation}\label{Eq:1}
d(\V_0) \subset (\C^2 \otimes \V_0) \oplus (\Lambda^2 \C^2 \otimes \V_1).
\end{equation}
Indeed if not then any element $v_0 \in \V_0$ not satisfying (\ref{Eq:1}) will generate
a projective $\widetilde{A}_{\Gamma,\zeta_{\C}}$-module $N_1$ and $M \simeq
M/N_1$ in $\UMod_2(\widetilde{A}_{\Gamma,\zeta_{\C}})$. Similarly,
we can assume that
\begin{equation}\label{Eq:2}
d(\V_1) \subset (\C^2 \otimes \V_1) \oplus (\Lambda^2 \C^2 \otimes \V_0)
\oplus \W.
\end{equation}
But then $d(\Lambda^2 \C^2 \otimes \V_1)=0$, and in particular,
the composition of maps $\V_1 \stackrel{d}{\to} \C^2 \otimes \V_1
\stackrel{d}{\to} \Lambda^2\C^2 \otimes \V_1$ endows $\Lambda\C^2 \otimes \V_1$
with a structure of $\widetilde{A}_{\Gamma,\zeta_{\C}}$-module.
By Lemma~\ref{L:4.61},
this implies that $\V_1=0$. 

\paragraph{}Let us fix some decomposition
  \begin{equation}\label{E:7}
\V=\bigoplus_a V_a \otimes \rho_a
\end{equation}
and set $V=\bigoplus_a V_a$. Let us split the map $d:\;\V \to \C^2 \otimes \V
\oplus \W$ as $d=d_0+d_1$ where $d_0: \V \to \C^2 \otimes \V$ and
$d_1: \V \to \W$. Then the maps $d_0$ and $d_1$ give rise, via
the identification
(\ref{E:7}), the fixed intertwiners (\ref{E:1.5}) and the map $u$, to elements
$B =\bigoplus_{h \in H} x_h \in E(V,V)$ and $j \in L(V,W)$ respectively.
Similarly,
the map $d:\; \C^2 \otimes \V \oplus \W \to \Lambda^2 \C^2 \otimes \V
\simeq \V$ gives rise to elements $C=\bigoplus_{h \in H} y_h \in E(V,V)$
and $i \in L(W,V)$. From the relation $dz=-zd$ for $u \in \C^2$ we deduce
that $y_h=\epsilon(h)x_h$ where $\epsilon:H \to \C^*$ is defined in
Section~8.3.
Similarly, from $d^2=\sum_a \zeta_{\C,a} c_a$ we deduce the relation 
$\mu(B,i,j)=\boldsymbol{\zeta}_{\C}$.

\paragraph{}Note that the assignment $M \to (B,i,j)$ depends on a choice of the
decomposition (\ref{E:7}), but that two such decompositions give rise to the
same element in $\mathcal{M}_{\zeta}(\mathbf{v},\mathbf{w})$.
Hence we have obtained in this way a well-defined map from the set of isomorphism
classes of objects in $\mathcal{N}_{\zeta}(\mathbf{w})$ to $\sqcup_{\mathbf{v}}
\mathcal{M}_{\zeta}(\mathbf{v},\mathbf{w})$. Conversely,
it is clear that any point $(B,i,j) \in \mathcal{M}_{\zeta}
(\mathbf{v},\mathbf{w})$ gives rise, via the above construction, to an
$\widetilde{A}_{\Gamma,\zeta_{\C}}$-module structure on the $A_\Gamma$-module
$M=\Lambda (\C^2
\otimes \V)[1] \oplus \W$. Moreover this module $M$ is equipped with a canonical
isomorphism $u: \W \stackrel{\sim}{\to} R(M)$.
  Thus $M \in \mathcal{N}_{\zeta}(\mathbf{w})$. This
map assigns distinct points in $\mathcal{M}_{\zeta}
(\mathbf{v},\mathbf{w})$ to nonisomorphic objects in $\mathcal{N}_{\zeta}(\mathbf{w})$
and the Theorem follows. \qed

\vspace{.1in}

\noindent
\textbf{Remark.} Theorem~\ref{P:4.61} is equivalent to Proposition~\ref{bij22}.

\section{Koszul duality and sheaves on $\mathbb{P}^2$}

\paragraph{9.1.} In this section we clarify the relation between the construction
of Section~8 and the moduli space of coherent
sheaves on some noncommutative deformations of the projective plane,
as studied in \cite{BGK}.

\paragraph{}Let $\UMod^0(\widetilde{A}_\Gamma)$ be the
full subcategory of $\UMod(\widetilde{A}_\Gamma)$ consisting of objects $V$ satisfying
\begin{equation}\label{E:5}
Ext^{l+n}_{\widetilde{A}_\Gamma}(\C,V[-l]) \neq 0 \Rightarrow n=0,
\end{equation}
where $\C$ denotes the trivial module.
The category $\UMod^0(A_\Gamma)$ is defined in a similar way.

\paragraph{}To any $\widetilde{A}_\Gamma$ module
$V =\bigoplus_{n\in \Z} V_n$ is associated a
complex
of $\Gamma$-equivariant coherent sheaves on $\mathbb{P}^2$
  $$\cdots \stackrel{d}{\to} L_i(V^*) \stackrel{d}{\to}
  L_{i+1}(V^*) \stackrel{d}{\to} \cdots$$
where $L_i(V^*)=(V^*)_i \otimes \mathcal{O}(i)$ and where the differential is
defined by
$(d\zeta)(x)=x\cdot \zeta(x)$ for $x \in \mathbb{P}^2$
and any section $\zeta \in \Gamma(L_i(V)).$ Here $V^*$ denotes the
$\widetilde{A}_{\Gamma}$-module dual to $V$.
A well-known theorem of Bernstein, Gelfan'd and Gelfan'd \cite{BGG} 
generalizing
the classical Koszul duality between $\C[x,y,z]$ and $\Lambda \C^3$  asserts
that the assignment $\Phi:V \to L_{\bullet}(V^*)$ induces an equivalence
between $\UMod(\widetilde{A}_\Gamma)^{op}$ and the
$\Gamma$-equivariant
derived category $D^b_\Gamma(Coh(\mathbb{P}^2))$
of coherent sheaves on $\mathbb{P}^2$. In particular under this equivalence
condition (\ref{E:5}) corresponds to
$H^i(L_{\bullet}(V^*))\neq 0 \Rightarrow i=0$ (see \cite{BGS}, Section 2.13),
and $\UMod^0(\widetilde{A}_\Gamma)^{op}$
is equivalent to the category
$Coh_{\Gamma}(\mathbb{P}^2)$ of $\Gamma$-equivariant coherent sheaves on
$\mathbb{P}^2$. Therefore, it is an abelian category.
Similar statements hold for $\UMod(A_\Gamma)$
and $D^b_\Gamma(\mathbb{P}^1)$.

\paragraph{}The functor $R$ restricts to a functor $\UMod^0(\widetilde{A}_\Gamma)
\to \UMod^0(A_\Gamma)$. It corresponds to the functor
of restriction 
$$D^b_\Gamma(Coh(\mathbb{P}^2)) \to
D^b_\Gamma(Coh(\mathbb{P}^1))$$
induced by the embedding $\mathbb{P}^1
\simeq \mathbb{P}((\C^2)^*) \hookrightarrow \mathbb{P}((\C^2 \oplus \C d)^*)
\simeq \mathbb{P}^2$. 

\paragraph{}Let us denote by $\Pi_\Gamma$ the preprojective algebra of the affine
quiver $(I,E)$
(see e.g \cite{Maffei}). It is well-known 
and easy to
check that $A(\Gamma)$ is Koszul and that its quadratic dual is 
$\Pi_\Gamma$
(see e.g \cite{HK}). Thus altogether we get the following diagram relating
various algebras :

\begin{equation*}
    \begin{CD}
      \Lambda\C^2 \rtimes \C[\Gamma] @>{Koszul\;duality.}>> \C[x,y] 
\rtimes \C[\Gamma]  \\
      @V{Morita\;eq.}VV      @VV{Morita\;eq.}V  \\
      A(\qQ) @>{Koszul\;duality}>> \Pi_\Gamma
    \end{CD}
\end{equation*}

Similarly, there is a diagram

\begin{equation}\label{Diagrr:1}
    \begin{CD}
    \Lambda \C^3 \rtimes \C[\G]  @>{Koszul\;duality.}>> \C[x,y,z] 
\rtimes \C[\Gamma]  \\
      @V{Morita\;eq.}VV      @VV{Morita\;eq.}V  \\
      \widetilde{A(\qQ)} @>{Koszul\;duality}>> \Pi_\Gamma[z]
    \end{CD}
\end{equation}

Our construction in Section~8 is based on a deformation of the left column of
diagram (\ref{Diagrr:1}). The corresponding right column consists of the
homogeneous coordinate ring of the noncommutative $\mathbb{P}^2$ studied
in \cite{BGK}, and the (graded) deformed preprojective algebra (see \cite{CBH}).

\paragraph{9.2.} In this paragraph we show how to recover the description
of Nakajima varieties as moduli space of torsion-free sheaves on $\mathbb{P}^2$ with 
fixed framing at $\infty$, using the representation theory of $\widetilde{A}_{\G}$.
This corresponds to case 1.3.i), i.e $\zeta_{\C}=0$ and $\zeta_{\R} \in (\mathbb{N}^+)^I$.

\vspace{.1in}

Let $\mathcal{T}$ be the full subcategory of
$\UMod^0(\widetilde{A}_\Gamma)$ consisting of modules
$T$ such that $\Phi(T)$
is a torsion sheaf on $\mathbb{P}^2$. We will say that an object $M$
of $\UMod^0_{\Z}
(\widetilde{A}_\Gamma)$ is \textit{torsion-free} if for any $T \in
\mathcal{T}$ we have
$\mathrm{Hom}_{\UMod^0(\widetilde{A}_\Gamma)}(M,T)=0$.

\begin{lem}\label{L:4.41}
Let $N$ be a graded $\widetilde{A}_\Gamma$-module such that
$N=\Lambda \C^2 \otimes \V$ as a $A_\Gamma$-module. If $\V \neq \{0\}$ then
$H^0(\Phi(N))$ is a nontrivial torsion sheaf.\end{lem}
\noindent
\textit{Proof.} This follows from a direct computation.\qed
\begin{lem}\label{L:4.42}
Let $N \in \UMod^0(\widetilde{A}_\Gamma)$ such that
$N_i =0$ for $i >0$ and such that $R(N) \simeq 0$. Then $N \simeq 0$.
\end{lem}
\noindent
\textit{Proof.} Suppose $N \not\simeq 0$. By Lemma~\ref{L:4.2},
$N$ decomposes as a
$A_\Gamma$-module as $N = \bigoplus_{i=r}^{-2} \Lambda \C^2 \otimes \V_i$
for some $\Gamma$-modules $\V_i$ and some $r \leq -2$. We can assume
\begin{equation}
d(\V_r)\subset \C^2 \otimes \V_r. \tag{a}
\end{equation}
Indeed, any $v_r \in \V_r$ not satisfying (a) generates a projective
$\widetilde{A}_\Gamma$-module $N'$ and $N \simeq N /N'$ in
$\UMod^0(\widetilde{A}_\Gamma)$. But then it follows from the
previous lemma that $H^{-r}(\Phi(N)) \neq 0$, in contradiction with the
assumption that $N \in \UMod^0(\widetilde{A}_\Gamma)$.\qed

\vspace{.1in}

As in Section~8.4 we us fix $\mathbf{w} \in \N^I$ and associate to it a 
$\G$-module $\mathbb{W}$. We will regard also $\W$ as a graded $A_\Gamma$-module, where
$\Lambda\C^2$ acts trivially, placed in degree $0$. Note that $\W$ is naturally
an object of $\UMod^0(A_\Gamma)$.
Let $\mathcal{N}(\mathbf{w})$ denote the set of pairs $(u,M)$ where
$M \in \UMod^0(\widetilde{A}_\Gamma)$ is torsion-free and
$u: \W \stackrel{\sim}{\to}
R(M)$ is an isomorphism.
We say two elements
$(u_1,M_1)$ and $(u_2,M_2)$ are isomorphic if there exists an isomorphism
$j:\;M_1 \to M_2$ such that $u_2=R(j) \circ u_1$.
\begin{theorem}\label{P:4.51}
There is a natural bijection between the set of isomorphism
classes of elements in $\mathcal{N}(\mathbf{w})$ and the set of points of
$\sqcup_\mathbf{v} \mathcal{M}_{\zeta}(\mathbf{v},
\mathbf{w})$.\end{theorem}
\noindent
\textit{Proof.} Let $(u,M) \in \mathcal{N}(\mathbf{w})$.
We first show the following
\begin{lem}\label{L:4.51}
There exists a $\Gamma$-module $\V$ and $M' \in
\UMod(\tilde{A}_\Gamma)$ such that $M \simeq M'$ and $M'=
\Lambda \C^2 \otimes \V [+1] \oplus \W$ as an $A_\Gamma$-module.\end{lem}
\noindent
\textit{Proof.} By Lemma~\ref{L:4.2}, there exists $\Gamma$-modules $\V_i$,
$i \in \Z$
such that $M \simeq \bigoplus_i \Lambda\C^2 \otimes \V_i[-i] \oplus \W$ as an
$A_\Gamma$-module. Observe that
\begin{equation}\label{Eqq:1}
d(\V_{-1}) \subset \Lambda^2 \C^2 \otimes \V_{-2} \oplus \C^2
\otimes \V_{-1} \oplus \W.
\end{equation}
Indeed if not then any $v_{-1} \in \V_{-1}$ such that (\ref{Eqq:1}) does not hold
generates a projective submodule $N$ of $M$, and $M \simeq M/N$ in
$\UMod^0(\tilde{A}_\Gamma)$.
Let $T$ be the $\widetilde{A}_\Gamma$-module obtained by restricting the
$\widetilde{A}_\Gamma$-action to $\bigoplus_{i \geq 0} \Lambda\C^2 \otimes
\V_i[-i]$.
Note that $T \in \mathcal{T}$. Indeed, we have $H^i(\Phi(T))=0$ for $i>0$ and
$H^i(\Phi(T))=H^i(\Phi(M))=0$ for $i <0$, so that
$T \in \UMod^0(\widetilde{A}_\Gamma)$, and $H^0(\Phi(T))$ is
torsion by
Lemma~\ref{L:4.41}. But $M$ is assumed to be torsion-free. This forces
$T \simeq 0$.\\
\hbox to1em{\hfill}A reasoning similar to (\ref{Eqq:1}) shows that
\begin{equation}
d(\V_{-2}) \subset \C^2 \otimes \V_{-2}
\oplus \Lambda^2 \C^2 \otimes \V_{-3}. \tag{b}
\end{equation}
Now observe that $N=\bigoplus_{i <-1} \Lambda \C^2 \otimes \V_i$ is in
$\UMod^0(\widetilde{A}_\Gamma)$ and that $R(N)\simeq 0$.
Thus $N \simeq 0$
by Lemma~\ref{L:4.42}, and the lemma follows.\qed

\paragraph{} Note that
by Lemma~\ref{L:4.42} again, the following holds :
\begin{equation}\label{E:6}
\mathrm{for\;all\;}N \subset M,\;R(N)=0 \Rightarrow N \simeq 0.
\end{equation}
and that condition (\ref{E:6}) is equivalent to the
stability condition of 1.3.i).

The rest of the proof of the Theorem now closely follows the proof of
Theorem~\ref{P:4.61}.\qed

%Moreover, this module $M$ is torsion free, in
%$\UMod^0(\widetilde{A}_\Gamma)$, and is equipped with a canonical
%isomorphism $u: \W \stackrel{\sim}{\to} R(M)$.

\paragraph{Remark.}
The above description of torsion-free sheaves on $\mathbb{P}^2$
with fixed framing at infinity is equivalent to the classical one in terms of
Beilinson's monads, (see e.g, \cite{Nak3}, Chap. 2) but its 
derivation doesn't use spectral sequences.

\paragraph{}From the above proof one easily deduces the following result.

\begin{corollary} The variety $\mathcal{M}_{\zeta}(\mathbf{v},\mathbf{w})$
is isomorphic to the set of all $\Lambda \C^2 \rtimes \C[\G]$-derivations 
(of degree one)
of the module $\Lambda \C^2 \otimes \mathbb{V} \oplus \mathbb{W}$ 
satisfying the
following condition : if $\mathbb{V}' \subset \mathbb{V}$ is a 
$\Gamma$-submodule
such that $\Lambda \C^2 \otimes \mathbb{V}'$ is $d$-stable then 
$\mathbb{V}'=0$.
\end{corollary}

\paragraph{9.3.} Denote by $\iota_a: \mathrm{End}(\rho_a) \to \C[\Gamma] 
\simeq
\bigoplus_a \mathrm{End}(\rho_a)$ the canonical embedding, and let
$\pi_a: \C[\G] \to \mathrm{End}(\rho_a)$ be the canonical projection.
We call 
$$m:\; \Lambda \C^2 \otimes \C[\G] \otimes \Lambda \C^2 \to \Lambda \C^2 
\otimes \C[\G]$$
the multiplication map and we define
$$\Delta:\;\Lambda \C^2 \otimes \C[\G] \to \Lambda \C^2 \otimes \C[\G]
\otimes \Lambda \C^2$$
to be its adjoint. Consider the following maps of 
$A_\Gamma$-bimodules :
$$d_1:\; \Lambda \C^2 \otimes \rho_a \otimes_\C \rho_a^* \otimes \Lambda \C^2
\simeq \Lambda \C^2 \otimes \mathrm{End}(\rho_a) \otimes \Lambda \C^2
\stackrel{m \circ (1 \otimes \iota_a \otimes 1)}{\longrightarrow}
\Lambda \C^2 \otimes \C[\G]$$
$$d_2:\; \Lambda \C^2 \otimes \C[\G] \stackrel{(1 \otimes \pi_a \otimes 1) 
\circ \Delta}{\longrightarrow} \Lambda \C^2 \otimes \mathrm{End}(\rho_a)
\otimes \Lambda \C^2.$$
As in Section~6, this gives rise, for any $x \neq 0$ to a duplex of 
$A_\G$-bimodules
$$\mathbf{C}_{a,x}:\; \qquad \stackrel{d_2}{\to} \Lambda \C^2 \otimes 
\mathrm{End}(\rho_a) \otimes \Lambda \C^2 \stackrel{d_1}{\to} \Lambda \C^2 
\otimes \C[\G] \stackrel{d_2}{\to}.$$
One checks that the $A_\Gamma$ bimodule duplex $\mathbf{C}_{a,x}$
corresponds to the $A(\qQ)$-bimodule duplex $C_{a,x}$ under the
equivalence $\mathbf{Mod}(A_\Gamma) \simeq \mathbf{Mod}(A
(\qQ))$. In particular, the collection of duplexes 
$\mathbf{C}_{a,x}$ satisfy the braid relations of Section~6.3, (in the
stable category of bimodule duplexes). Thus, as in Section~7, tensoring by 
$\mathbf{C}_{a,x}$ for $a \in I$ and generic $x$ defines an
action of the Weyl group $\mathbf{W}$ on the set of objects of
$\UMod_2(A_{\Gamma,\zeta_{\C}})$ for all generic
$\zeta_{\C}$. 

\paragraph{} In other words, (for generic $\zeta_{\C})$ and
$a \in I$ we have a functor
$$\mathcal{R}_a: \UMod_2(\widetilde{A}_{\Gamma,\zeta_{\C}})
\to \UMod_2(\widetilde{A}_{\Gamma, s_a(\zeta_{\C})}),$$ 
and the collection of such functors satisfy the braid relations.
Moreover, there is a canonical natural transformation 
$R \circ \mathcal{R}_a \to R$, and for any given fixed $\mathbf{w} \in 
\N^I$, $\mathcal{R}_a$ acts on
the set of objects of $\mathcal{N}_{\zeta}(\mathbf{w})$. 
The following Proposition is a consequence
of Theorem~\ref{T:2}.

\begin{prop} The action of $\mathcal{R}_a$ on 
$\mathcal{N}_{\zeta}(\mathbf{w})$ coincides 
with the action of $\kappa_a$ under the 
identification $\mathcal{N}_{\zeta}(\mathbf{w}) \simeq
\sqcup_{\mathbf{v}} \mathcal{M}_{\zeta}
(\mathbf{v},\mathbf{w}).$ 
\end{prop}

\small{}
\vspace{4mm}
\noindent
I. Frenkel, Dept. of Mathematics, Yale University, 10 Hillhouse Ave.
New Haven CT 06520-8283 USA; \vspace{0.1in}\\  
M. Khovanov, Dept. of Mathematics, University of California, One Shields Avenue
Davis, CA 95616-8633 USA, \texttt{mikhail@math.ucdavis.edu}\vspace{0.1in}\\
O. Schiffmann, Dept. of Mathematics, Yale University, 10 Hillhouse Ave.
New Haven CT 06520-8283 USA; and DMA Ens Paris, 45 rue d'Ulm, 75005 Paris
FRANCE, \texttt{schiffma@dma.ens.fr}.


\begin{thebibliography}{99}
\bibitem[BGK]{BGK}
Baranovsky V., Ginzburg V., Kuznetsov A., \emph{Quiver varieties and a
noncommutative $\mathbb{P}^2$}, Compositio Math. \textbf{134},  no. 3, 283--318, (2002).
\bibitem[BGS]{BGS}
Beilinson A., Ginzburg V., Soergel W., \emph{Koszul duality patterns in 
representation theory}, J. Amer. Math. Soc. \textbf{9}, no. 2, 473--527 (1996).
\bibitem[BGG]{BGG}
Bernstein I.N., Gelfand I.M., Gelfand S.I., \emph{ Algebraic vector 
bundles on $P\sp{n}$ and problems of linear algebra.} (Russian)
Funktsional. Anal. i Prilozhen. \textbf{12}, no.3, 66--67, (1978).
\bibitem[Bo]{Bo} Borcherds R., \emph{Generalized Kac-Moody algebras,} Journal of 
Algebra {\bf 115} 2 (1988), 501--512. 
\bibitem[CH]{CBH}
Crawley-Boevey W., Holland M., \emph{Noncommutative deformations of Kleinian singularities},
Duke Math. J. \textbf{92},  no. 3 605--635 (1998). 
\bibitem[H]{Happel}
Happel, D., \emph{Triangulated categories in the representation theory of
finite-dimensional algebras.},
London Mathematical Society Lecture Note Series, \textbf{119}, (1988).
\bibitem[HK]{HK}
Huerfano, R.S., Khovanov, M., \emph{A category for the adjoint representation.}
J. Algebra \textbf{246}, no. 2, 514--542, (2001), arXiv math.QA/0002060.  
%\bibitem[K]{Kac} 
%Kac, V.G., \emph{Infinite dimensional Lie algebras,} Cambridge U Press, 1990.  
\bibitem[KS]{KS}
Khovanov M., Seidel P., \emph{Quivers, Floer cohomology, and braid group 
actions.}, J. Amer. Math. Soc. \textbf{15}, no. 1, 203--271, (2002), 
 arXiv math.QA/0006056.  
\bibitem[K]{Ku} Kumar S., \emph{Kac-Moody groups, their flag 
varieties and representation theory.} Birkha\"user, Boston, 2002.   
\bibitem[L]{LusWeyl}
Lusztig G., \emph{Quiver varieties and Weyl group actions.},
  Ann. Inst. Fourier \textbf{50}, no. 2, 461--489, (2000).
\bibitem[M]{Maffei}
Maffei A., \emph{A remark on quiver varieties and Weyl groups}, arXiv 
 math.AG/0003159.  
\bibitem[MS]{MS} McDuff D., Salamon D., 
 \emph{J-holomorphic curves and quantum cohomology}, University Lecture Series 6, 
 AMS, 1994.  
\bibitem[N1]{Nak1}
Nakajima H., \emph{Instantons on ALE spaces, quiver varieties, and Kac-Moody algebras}, 
Duke Math. J. \textbf{76}, no. 2, 365--416, (1994).
\bibitem[N2]{Nak2}
Nakajima H., \emph{Reflection functors for quiver varieties and Weyl group
actions}, Math. Ann., to appear.
\bibitem[N3]{Nak3}
Nakajima H., \emph{Lectures on Hilbert scheme of points on 
surfaces}, Uni. Lect.
Ser. AMS,  n. \textbf{18}, (1999).
\bibitem[PX]{PX} Peng, L., Xiao, J., 
\emph{Triangulated categories and Kac-Moody algebras},  Invent. Math.  \textbf{140},  
no. 3, 563--603, (2000) 
\bibitem[RZ]{RZ}
Rouquier R., Zimmermann A., \emph{Picard groups for derived module 
categories}, Proc. London Math. Soc. (3)  \textbf{87},  no. 1, 197--225, (2003).
\bibitem[ST]{ST}
Seidel P., Thomas R., \emph{Braid group actions on derived categories of 
coherent sheaves.}, Duke Math. J. \textbf{108}, no. 1, 37--108, (2001), 
 arXiv math.AG/0001043.  
\end{thebibliography}
\end{document}